\documentclass[12pt]{article}
\usepackage{xcolor}
\usepackage[colorlinks]{hyperref}
\usepackage{hyperref}
\hypersetup{
	linkcolor= blue
	,citecolor= green
	,filecolor= cyan
	,urlcolor= magenta
	,menucolor= red
	,runcolor= cyan
}
\usepackage{float}
\usepackage{authblk}
\usepackage[margin=1in]{geometry}
\usepackage{graphicx}
\usepackage{subcaption}% Required for inserting images
\usepackage{amsmath}
\usepackage{accents}
\usepackage{amssymb}
\usepackage{amsthm}
\usepackage{booktabs}
\usepackage{bm}
\usepackage{cleveref}
\usepackage[dvipsnames]{xcolor}
\usepackage{tikz}
\usepackage{multirow}
\usetikzlibrary{shapes,arrows.meta,positioning}

\newcommand{\eqautoref}[1]{\hyperref[#1]{Eq.~(\ref*{#1})}}

\title{Physics-Informed Functional Link Constrained Framework with Domain Mapping for Solving Bending Analysis of an Exponentially Loaded Perforated Beam}
\author[1]{Iswari Sahu\thanks{\textcolor{magenta}{\texttt{sahuiswari920@gmail.com}}}}
\author[1]{Ramanath Garai\thanks{\textcolor{magenta}{\texttt{ramanathgarai12@gmail.com}}}}
\author[1]{S. Chakraverty\thanks{\textcolor{magenta}{\texttt{Corresponding author:sne\_chak@yahoo.com}}}}
\affil[1]{Department of Mathematics, National Institute of Technology Rourkela, India}
\date{}
\begin{document}
\maketitle
	\begin{abstract}
 This article presents a novel and comprehensive approach for analyzing bending behavior of the tapered perforated beam under an exponential load. The governing differential equation includes important factors like filling ratio ($\alpha$), number of rows of holes ($N$), tapering parameters ($\phi$ and $\psi$), and exponential loading parameter ($\gamma$), providing a realistic and flexible representation of perforated beam configuration.
Main goal of this work is to see how well the Domain mapped physics-informed Functional link Theory of Functional Connection (DFL-TFC) method analyses bending response of perforated beam with square holes under exponential loading.
 For comparison purposes, a corresponding PINN-based formulation is developed. Outcomes clearly show that the proposed DFL-TFC framework gives better results, including faster convergence, reduced computational cost, and improved solution accuracy when compared to the PINN approach. These findings highlight effectiveness and potential of DFL-TFC method for solving complex engineering problems governed by differential equations. 
 Within this framework, hidden layer is replaced by a functional expansion block that enriches  input representation via orthogonal polynomial basis functions, and the domain of DE mapped to corresponding domain of orthogonal polynomials.
 A Constrained Expression (CE), constructed through the Theory of Functional Connections (TFC) using boundary conditions, ensures that constraints are exactly satisfied. In CE, free function is represented using a Functional Link Neural Network (FLNN), which learns to solve resulting unconstrained optimization problem.  The obtained results are further validated through the Galerkin and PINN solutions. 
\end{abstract}
\textbf{Keywords:} DEs;  DFL-TFC; Perforation; Bending deflection; Exponential Load; Taper beam

\section{Introduction}
Beams play a crucial role in engineering systems because they efficiently transfer and resist bending loads. Their length is typically much larger than the cross-sectional dimensions, which allows them to act as primary load-carrying members in many structures. In modern engineering designs, beams are sometimes manufactured with periodically arranged openings, forming what are known as perforated beams. These holes are often described through a filling ratio that quantifies the portion of solid material within each periodic segment. Presence of such perforations can significantly reduce structural weight and material usage while simultaneously altering stiffness and mass distribution of the beam. As a result, understanding the mechanical behavior of perforated beams has become an important topic in structural analysis, particularly when their response is governed by complex differential equations that may involve nonlinear, coupled, or mixed-order terms.

To analyze these complex beam problems, conventional numerical methods such as the finite element method (FEM) \cite{zienkiewicz1977finite} and the finite difference method (FDM) \cite{godunov1959finite} are commonly employed. Although these approaches are widely used, they may involve high computational costs and convergence difficulties when applied to high-order differential equations or complex structures such as perforated beams. In recent years, scientific machine learning has emerged as a promising alternative by integrating machine learning techniques with scientific computing to approximate PDE solutions. However, purely data-driven approaches usually require large datasets and can be sensitive to noise. To address these limitations, recently, machine learning approaches, particularly Physics-Informed Neural Network (PINN) \cite{raissi2019physics} has been developed to solve a range of differential equations (DEs). PINN incorporates governing equations into the loss function and eliminate the need for labeled datasets. However, they enforce boundary conditions using penalty terms, which can lead to slow convergence, sensitivity to hyperparameters, and training instability. Several review studies \cite{cuomo2022scientific,blechschmidt2021three} provide a comprehensive overview of recent advances in scientific machine learning and highlight its applications to real-world engineering problems.

Luschi and Pieri \cite{luschi2012simple,luschi2014analytical} developed analytical models to investigate the resonance frequencies of perforated beams, examining beams with periodic rectangular and square perforations and deriving closed-form expressions for their equivalent bending and shear stiffness. Accordingly, Fallah et al. \cite{fallah2024physics} applied physics-informed neural networks (PINNs) to study the bending and free vibration of three-dimensional functionally graded porous beams resting on an elastic foundation, with material properties varying continuously in three spatial directions. The forward and inverse bending response of nanobeams resting on a three-parameter nonlinear elastic foundation, within the framework of nonlocal elasticity theory, was investigated using physics-informed neural networks (PINNs) by Kianian et al. \cite{kianian2025pinn}.
Chen and  Liang \cite{chen2023second} presented a PINN approach for second-order beam-column analysis, delivering a mesh-free, stable substitute for standard finite element method (FEM) in handling large deflections. Bazmara and Mianroodi \cite{bazmara2023application} solved the nonlinear large-deformation buckling problem of 3D FG porous beams on a Winkler-Pasternak foundation using PINN. 
Teloli et al. \cite{teloli2025physics} solved inverse structural problems of the fourth-order Euler-Bernoulli beam equation using PINN. This framework predict displacements and find structural factors like damping and elastic modulus. Kapoor and Wang \cite{kapoor2024transfer} developed a causality-respecting PINN combined with transfer learning to simulate Euler-Bernoulli and Timoshenko beams on a Winkler foundation. 
Although PINN have shown strong potential, they face several challenges in certain cases \cite{wang2022and, krishnapriyan2021characterizing}. They may struggle with optimization during training, experience slow convergence, and incur high computational costs, particularly for complex or high-dimensional problems.

In this work, we analyze the bending behavior of perforated beam structures in the presence of exponential load using the propose DFL-TFC method. The beam is considered to have a non-uniform cross-section, with its height varying along the axial direction, allowing for a comprehensive analysis of tapered structural configurations under complex loading conditions. DFL-TFC is a hybrid framework integrating FLNN and TFC, where the DE domain is mapped onto an orthogonal polynomial domain. In the FLNN framework, the input features are then expanded into a higher-dimensional feature space using orthogonal polynomials, rather than employing hidden layers. Owing to its minimal number of trainable parameters, the FLNN is computationally efficient. Meanwhile, TFC constructs a constrained expression (CE) that inherently satisfies the initial and boundary conditions (ICs and BCs), embedding these constraints directly into the solution structure through CE, the CE removes the necessity of adding penalty terms in the loss function. It delivers accurate results using fewer trainable parameter than PINNs.

The rest of the paper is structured as follows: 
 Section \ref{sec.Perforated Beam} presents theoretical modeling of the tapered perforated beam, including the governing equations and underlying assumptions.
Section \ref{sec. problem dev} describes the formulation of the problem, outlining the domain, boundary conditions, and relevant parameters.
Section \ref{methodology} describes the methodology of PINN, DFL-TFC, and Galerkin procedure.
Section \ref{sec. results} presents the validation of the DFL-TFC method and gives results along with a detailed discussion.
Finally, Section \ref{sec. conclusion} concludes the paper and highlights potential directions for future research.

\section{Theoretical Modeling of Perforated Beam}
\label{sec.Perforated Beam}
\begin{figure}[h]
		\centering
		\includegraphics[height=2.7in]{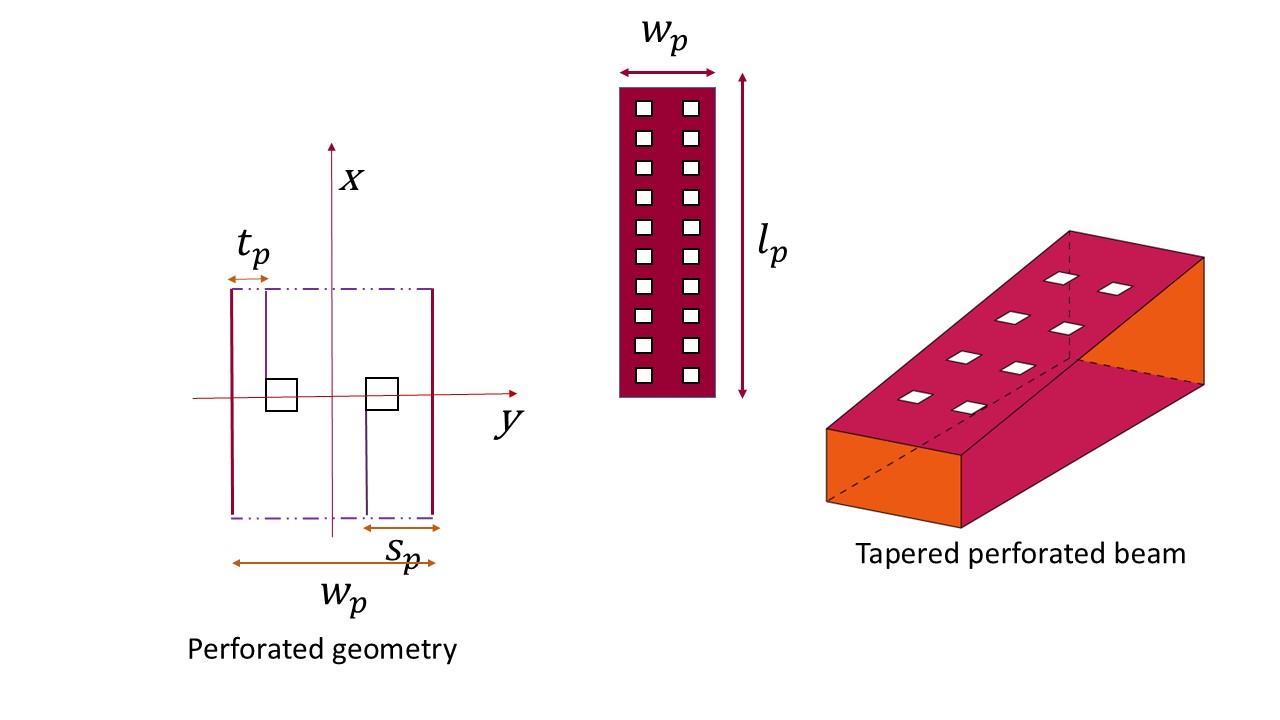}
		\caption{ Geometric model of a tapered perforated beam }
		\label{p2perfo1}
	\end{figure}

A perforated beam is an engineered structural member in which material is selectively removed through a planned pattern of perforations to improve its mechanical performance. These perforations, which may be circular, rectangular, or of other geometries, influence the stiffness, mass distribution, and dynamic response of the beam. In the present configuration, the beam incorporates square-shaped perforations arranged periodically in a grid pattern.\\ 
	The geometric configuration of a perforated beam is illustrated in Fig.~\ref{p2perfo1}. The beam has an overall length $l_p$, width $w_p$, and thickness $h_p$. Square perforations are periodically distributed along the beam, defined by a spatial period $s_p$ and period length $t_p$. Accordingly, the side length of each square hole is $s_p - t_p$, and $N $ represents number of holes along the length of beam. The ratio of period length $t_p$ to spatial period $s_p$ is called the filling ratio of perforated beam, following the definitions reported in the literature \cite{luschi2014analytical,eltaher2018}.
	\begin{equation}
		\alpha = \frac{t_p}{s_p},  \quad 0 \leq \alpha \leq 1 .
	\end{equation}
	When $\alpha = 1$, the beam converted into a fully solid configuration with no material removal. For intermediate values $0 < \alpha < 1$, the beam contains perforations and is therefore partially filled. In the limiting case $\alpha = 0$, beam represents a completely perforated (idealized) structure. Hence, the filling ratio is expressed as
	$$
	\alpha =
	\begin{cases}
		0 & \text{Completely Perforated (limiting case)}, \\
		(0,1) & \text{Partially Filled}, \\
		1 & \text{Fully Solid}.
	\end{cases}
	$$

\subsection{Geometric Tapering of Perforated Beams}
In the present study, perforated beam is assumed to have a non-uniform cross-section, with the beam height varying smoothly along axial direction. The height variation is described by a quadratic function as
\begin{equation}
h(x) = h_0 \left( 1 + \phi \frac{x}{l_p} + \psi \left(\frac{x}{l_p}\right)^2 \right),
\label{eq:height_variation}
\end{equation}
where $h_0 $ denotes height of the beam at the left end ( $x = 0 $), $ l_p $ is total length of the beam, and \( \phi \) and \( \psi \) are tapering parameters governing linear and quadratic variations of the height, respectively.

This quadratic tapering profile ensures a continuous modification of the beam geometry, resulting in a gradual redistribution of bending stiffness and mass along the beam axis. Gradual changes in shape avoid sudden stiffness variations, which lowers stress buildup and leads to better static deflection.

\noindent Modified bending stiffness of the tapered perforated beam, accounting for effects of periodic perforations and tapered geometry, can be expressed as \cite{luschi2014analytical,akpinar2025size}:
	\begin{equation} \label{p2eq2}
		[EI]_{\text{eq}} = [E I(x)]_{s}
		\left[ \frac{\alpha (N+1)(N^2 + 2N + \alpha^2)}
		{\left(1 - \alpha^2 + \alpha^3\right)N^3 + 3\alpha N^2 + \left(3 + 2\alpha - 3\alpha^2 + \alpha^3\right)\alpha^2N + \alpha^3} \right],
	\end{equation}
where $E$ is the Young’s modulus and $I(x)$ is the second moment of area of the tapered beam. Here, $[EI]_{\mathrm{eq}}$ represents modified bending stiffness of the axially tapered perforated beam, accounting for the influence of periodic perforations, while $[EI]_s$ denotes the bending stiffness of corresponding solid tapered beam.\\
 Modified mass per unit length of the tapered perforated beam, relative to  corresponding fully solid beam, together with associated modified moment of inertia per unit length, can be formulated as \cite{alnujaie2023influence,abdelrahman2021},
	
	\begin{equation}
		[\rho A]_{eq} = [\rho A(x)]_s   \left[ \frac{[1 - N(\alpha - 2)]\alpha}{N + \alpha} \right],
	\end{equation}
	\begin{equation}
		[\rho I]_{\text{eq}}= [\rho I(x)]_s  \left[ \frac{\alpha \left[ (2 - \alpha)N^3 + 3N^2 - 2(\alpha - 3)(\alpha^2 - \alpha + 1)N + \alpha^2 + 1 \right]}{(N + \alpha)^3} \right].
	\end{equation}
Here, $\rho$ denotes the material density of beam. The term $[\rho A]_{\mathrm{eq}}$ represents the modified mass per unit length of the tapered perforated beam, while $[\rho I]_s$ denotes the moment of inertia per unit length of the corresponding solid tapered beam.

\section{Problem Development} \label{sec. problem dev}                    
The governing differential equation describing bending behavior of a tapered perforated beam resting on an elastic foundation and subjected to an exponential load can be expressed as \cite{kafkas2025size}:
\begin{equation} \label{coleq17}
\frac{d^2}{dx^2} \left[ [EI]_{eq} \frac{d^2 W_p(x)}{dx^2} \right] =  q(x)
+ K_{p} \frac{d^2 W_p(x)}{dx^2} ,
\end{equation}     
where $[EI]_{eq}$ denotes modified bending stiffness of the axially tapered perforated beam, $K_p$ is Pasternak foundation variable, and $q(x)$ represents external exponential load acting on the beam.

Here, the exponential load is taken in following forms:
$$q(x)= q_0 e^{\gamma \frac{x}{l_p}} .$$

After substituting the above external load into Eq. (\ref{coleq17}), it yields:
\begin{equation} \label{coleq19}
\frac{d^2}{dx^2} \left[ [EI]_{eq} \frac{d^2 W_p(x)}{dx^2} \right]  =  q_0 e^{\gamma \frac{x}{l_p}}
+ K_{p} \frac{d^2 W_p(x)}{dx^2}   .
\end{equation}

Introducing non-dimensional variables removes dependence on physical units and characteristic scales, resulting in a generalized formulation. Accordingly, the following non-dimensional parameters are defined for subsequent analysis.
$$X = \frac{x}{l_p}, \ \  \overline {W_p}(X)=\frac{W_p(x)}{l_p}, \ \ q(X) = \frac{q(x)}{l_p}$$

After using non-dimensional terms in Eq. (\ref{coleq19}), it yields:
\begin{equation}
    \frac{ 1}{l_p^2} \frac{d^2}{dX^2} \left[ \frac{[EI]_{eq}}{l_p} \frac{d^2 \overline W_p(X)}{dX^2} \right] = l_p q(X) + \frac{K_p}{l_p^2} \frac{d^2 \overline W_p(X)}{dX^2} .
\end{equation}

After simplification, it gives:
\begin{equation} \label{coleq21}
    \frac{d^2}{dX^2} \left[ \widetilde{[EI]}_{eq} \frac{d^2 \overline W_p(X)}{dX^2} \right] =  \widetilde{ q_0} e^{\gamma X}+ \widetilde{K_p} \frac{d^2 \overline W_p(X)}{dX^2},
\end{equation}
where $\widetilde{ q_0} = \frac{q_0 l_p^4}{E I} ,$
$\widetilde{K_p} = \frac{K_p l_p^2}{E I} ,$
$$ \widetilde{[EI]}_{eq} = (1+\phi X +\psi X^2)^3 \left[ \frac{\alpha (N+1)(N^2 + 2N + \alpha^2)} {\left(1 - \alpha^2 + \alpha^3\right)N^3 + 3\alpha N^2 + \left(3 + 2\alpha - 3\alpha^2 + \alpha^3\right)\alpha^2N + \alpha^3} \right] .$$

\subsection{Imposed Boundary Conditions}
Boundary conditions at the two ends of  tapered perforated beam, corresponding to $X=0$ and $X=1$, are specified in their most general linear form as follows \cite{garai2026effect}.
\begin{equation}
		c_1 \, \overline{W_p} + c_2 \, \frac{d\overline{W_p}}{dX} 
		+ c_3 \, \left[ [EI]_{\mathrm{eq}} \, \frac{d^2 \overline{W_p}}{dX^2} \right]
		+ c_4 \, \frac{d}{dX} \left[ [EI]_{\mathrm{eq}} \, \frac{d^2 \overline{W_p}}{dX^2} \right] = 0.
	\end{equation}
	\begin{equation}
		d_1 \, \overline{W_p} + d_2 \, \frac{d\overline{W_p}}{dX} 
		+ d_3 \,  \left[ [EI]_{\mathrm{eq}} \, \frac{d^2 \overline{W_p}}{dX^2} \right]
		+ d_4 \, \frac{d}{dX} \left[ [EI]_{\mathrm{eq}} \, \frac{d^2 \overline{W_p}}{dX^2} \right] = 0.
	\end{equation}

The constants $c_i \ \text{and}  \ d_i, \ i = 1,2,3,4 $  define boundary conditions at the ends of  tapered beam. Their values depend on the type of support used.

\noindent Simply Supported End:
\[
\begin{aligned}
	c_1 &= 1, \quad c_2 = 0, \quad c_3 = 0, \quad c_4 = 0 
	&&\Rightarrow \overline{W_p} = 0,\\
	d_1 &= 0, \quad d_2 = 0, \quad d_3 = 1, \quad d_4 = 0 
	&&\Rightarrow [EI]_{\mathrm{eq}}\,\overline{W_p}'' = 0.
\end{aligned}
\]

\noindent Clamped End:
\[
\begin{aligned}
	c_1 &= 1, \quad c_2 = 0, \quad c_3 = 0, \quad c_4 = 0 
	&&\Rightarrow \overline{W_p} = 0,\\
	d_1 &= 0, \quad d_2 = 1, \quad d_3 = 0, \quad d_4 = 0 
	&&\Rightarrow \overline{W_p}' = 0.
\end{aligned}
\]

\section{Methodology} \label{methodology}
This section describes the methodologies of Physics Informed Neural Network, Galerkin Procedure, and Functional Link Theory of Functional Connection Frameworks with Domain Mapping.

\subsection{Physics Informed Neural Network} \label{sec. PINN}
The governing differential equation (\ref{coleq21}) is written in the form:
\begin{equation*} 
    \frac{d^2}{dX^2} \left[ \widetilde{[EI]}_{eq} \frac{d^2 {\overline W_p}(X)}{dX^2} \right] =  \widetilde{ q_0} e^{\gamma X}+ \widetilde{K_p} \frac{d^2 {\overline W_p}(X)}{dX^2}.
\end{equation*}

In physics informed neural network (PINN) \cite{raissi2019physics,kumar2023physics}, input variable X is fed into the network and processes them through 3 hidden layers, each hidden layer contain 5 neuron equipped with tanh activation function. The information is passes through these layers, and the final layer produces the predicted output $\widehat{\overline{W}_p}(X)$.  In this study, networks with one, two, and three hidden layers are examined in the loss analysis Section \ref{loss analysis}. Based on the loss analysis, the network with three hidden layers provides most optimized results; therefore,  three hidden layer is taken here.

Residual of the DE is given as: 
\begin{equation}\label{residual}
 \mathcal{R}(X) = \frac{d^2}{dX^2} \left[\widetilde{[EI]}_{eq} \frac{d^2 \widehat{\overline W_p}(X)}{dX^2} \right] -  \widetilde{ q_0} e^{\gamma X}- \widetilde{K_p} \frac{d^2 \widehat{\overline W_p}(X)}{dX^2}  .
\end{equation}

The total loss function is given below:
\begin{equation}\label{eq. PINN loss}
    \mathcal{L}_{\text{Total}} =  \mathcal{L}_{\text{DE}} + \mathcal{L}_{\text{BC}} ,
\end{equation}
where
\begin{align*}
\mathcal{L}_{\mathrm{DE}} =
\frac{1}{N_f}
\sum_{i=1}^{N_f}
\left|
\frac{d^2}{dX^2} \left[\widetilde{[EI]}_{eq}  \frac{d^2 \widehat{\overline W_p}(X)}{dX^2} \right] -  \widetilde{ q_0} e^{\gamma X}- \widetilde{K_p} \frac{d^2 \widehat{\overline W_p}(X)}{dX^2} 
\right|^2,
\end{align*}
where $N_f$ is the number of collocation points.

Loss for the S-S boundary conditions is given by:
\begin{align*}
\mathcal{L}_{\text{BCs}} 
= \frac{1}{4} 
\left(\left| \widehat{\overline{W}_p}(0) \right|^2 + \left| \widehat{\overline{W}_p}(1) \right|^2 + \left| \widehat{\overline{W}''_p}(0) \right|^2 + \left| \widehat{\overline{W}''_p}(1) \right|^2\right) .
\end{align*}

 For C–S boundary conditions, loss is given by:
\begin{align*}
\mathcal{L}_{\text{BCs}} 
= \frac{1}{4} 
\left(\left| \widehat{\overline{W}_p}(0) \right|^2 + \left| \widehat{\overline{W}_p}(1) \right|^2 + \left| \widehat{\overline{W}'_p}(0) \right|^2 + \left| \widehat{\overline{W}''_p}(1) \right|^2\right) .
\end{align*}

Training is carried out using 100 points distributed over the domain [0,1].
Automatic differentiation (AD) \cite{paszke2017automatic}is used to compute the gradients of the neural network with respect to its trainable parameters.
 The trainable parameters, containing all weights and biases of network, is iteratively optimized by minimizing the total loss function defined in \eqautoref{eq. PINN loss}. This optimization is performed using gradient-based algorithms such as L-BFGS \cite{kingma2014adam, liu1989limited} with 50 optimization steps, where each step includes 50 internal iterations. After training process is completed, neural network offers an approximate solution to the governing differential \eqautoref{coleq21} and yields reliable predictions in the domain.

\begin{figure}[H]
    \centering
    \includegraphics[width=0.993\linewidth]{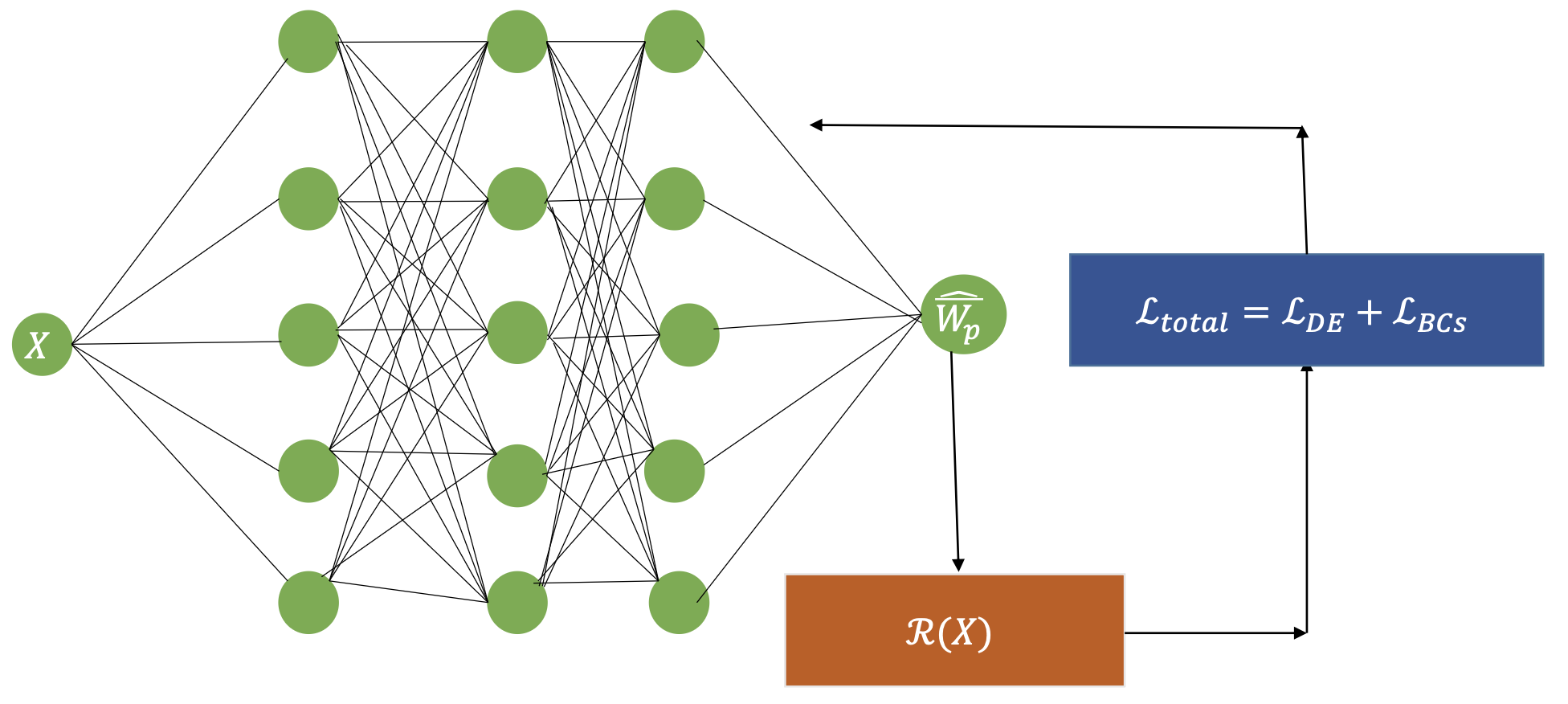}
    \caption{Schematic diagram of PINN}
    \label{PINN}
\end{figure}

\subsection{Galerkin Approximation Procedure} \label{galerkin_sec}
The transverse displacement of the perforated tapered beam, denoted by $ \overline{W}_p(X) $, is approximated using a finite series of basis functions as \cite{bhat1985natural,garai2026effect}:
\begin{equation} \label{coleq24}
\overline{W}p(X) = \sum_{i=1}^{n} \eta_i \boldsymbol{\Theta}_i(X),
\end{equation}
where $n$ is number of basis functions, $\eta_i$ are unknown coefficients, and $\boldsymbol{\Theta}_i(X)$ represent orthonormal polynomial functions, which comes from simple polynomial $X^{i-1}, \ i = 1,2,..n$ using the Gram-Schmidt process.

By applying one of the boundary conditions (S-S or C-S) into Eq. (\ref{coleq24}), the updated displacement function is written as:
\begin{equation} \label{coleq25}
\overline{W}p(X) = \sum_{i=5}^{n} \eta_i \boldsymbol{\overline \Theta}_i(X),
\end{equation}
where $\boldsymbol{\overline \Theta}_i(X)$ is the modified orthonormal polynomial.

Residual associated with the governing equation (\ref{coleq21}) is defined as:
\begin{equation}
    \mathfrak{R}(\overline{W}_p) = \frac{d^2}{dX^2} \left[ \widetilde{[EI]}_{eq}   \frac{d^2 \overline W_p(X)}{dX^2} \right] - \widetilde{ q_0} e^{\gamma X}  -\widetilde{K_p} \frac{d^2 \overline W_p(X)}{dX^2}.
\end{equation}
	
The weight function is taken as:
\begin{equation}
\mathfrak{w}_i = \boldsymbol{\overline {\Theta}}_i(X), \quad i = 5, 6, \ldots, n.
\end{equation}
Taking inner product of the residual with weight functions by the following expression:
\begin{equation}
\langle \mathfrak{R}(\overline{W}_p), \mathfrak{w}_i  \rangle = 0, \quad i=5,6,\ldots,n,
\end{equation}
This leads to a system of $n-4$ equations. Solving these $n-4$ governing equations yields the bending deflection $\overline{W}_p(X)$. The result is then multiplied by a factor of $100$ and written as:
$$\widetilde{W_p}(X) = 100 \times \overline W_p(X), \  \ X \in [0,1].$$

\subsection{Functional Link-Theory of Functional Connection Framework with Domain Mapping}\label{sec-FLTFC}
 Theory of functional connection \cite{mortari2017theory, leake2020multivariate, leake2020deep} provides a structured framework for building functions that automatically meet given constraints. This framework combines a free function $h(x)$, projection functionals $\boldsymbol{\rho}(x, h(x))$, and switching functions $\Phi_{x_i}$.
This methodology is effective for developing constrained formulations in both one-dimensional and multidimensional domains.

The general form of CE is define as \cite{leake2020multivariate, leake2020deep,sahu2026physics}
\begin{equation}\label{eq. CE formula}
    \hat{f}(\mathbf{x}, h(\mathbf{x})) = h(\mathbf{x}) +  \mathcal{M} \left( \boldsymbol{\rho}(\mathbf{x}, h(\mathbf{x})) \right)\, \Phi_{x_1}\, \Phi_{x_2} \cdots \Phi_{x_n},  
\end{equation}
where \(\mathcal{M}\) represents an \(n\)-dimensional tensor, a multidimensional array with components composed of constraint and associated projection functionals and   \(\mathbf{x} = (x_1, x_2, \dots, x_n)^\top \in \mathbb{R}^n\) is the input vector. 

To make the CE \cite{ leake2020multivariate, leake2020deep,sahu2026physics}, follow these steps:
\begin{itemize}
    
    \item If a single index varies in the \(k^{\text{th}}\) dimension while all other indices are fixed at 1, the associated tensor elements are expressed as 
    \begin{equation}\label{eq. m_11}
        \mathcal{M}_{1\cdots x_k \cdots 1} = 
    \left\{0, {}^{x_k}\rho_{1},\cdots,{}^{x_k}\rho_{l_{x_k}}\right\}.
    \end{equation}

    \item A tensor element where more than one index is not equal to 1 is define as
    \begin{equation}\label{eq. m_ij}
        \mathcal{M}_{i_{x_1}i_{x_2}...i_{x_n}} = (-1)^{n+1}\,{}^{x_1}\mathfrak{C}^{i_{x_1}-1}\bigg[{}^{x_2}\mathfrak{C}^{i_{x_2}-1}\Big[\cdots\left[{}^{x_n}\rho_{i_{x_n}-1}\right]\cdots\Big]\bigg].
    \end{equation}

    \item For each variable \(x_k\), the vector of switching functions of the constraints \(l_{x_k}\)  is given by
    \begin{equation}
        \Phi_{x_k} = \begin{bmatrix}
        1 \\ {}^{x_k}\phi_1 \\ \vdots \\ {}^{x_k}\phi_{l_{x_k}}
    \end{bmatrix} \in \mathbb{R}^{l_{x_k} + 1}
    \end{equation}

    \item  support functions \(l_{x_k}\) are linearly combined to provide each switching function \({}^{x_k}\phi_j\)
    \begin{equation}\label{eq. switching fn}
        {}^{x_k}\phi_j(x_k) = \sum_{i=1}^{l_{x_k}} \alpha_{ij}\, \mathtt{s}_i(x_k),
    \end{equation}
    where \(\mathtt{s}_i(x_k) = x_k^{i-1}\) are the monomial support functions.

    \item The support functions ${\mathtt{s}_i(x_k)}$ are chosen to maintain linear independence across the domain of $x_k$. This property guarantees that the associated coefficient matrix remains invertible.
 The coefficient matrix takes the form
    \begin{equation}\label{eq. coefficient matrix}
        \left[ \alpha_{ij} \right]_{l_{x_k} \times l_{x_k}} = \left[ {}^{x_k}\mathfrak{C}^i [\mathtt{s}_j(x_k)] \right]^{-1}.
    \end{equation}
    
  To determine the coefficients $\alpha_{ij}$, each switching function $\phi_j(x)$ must satisfy the following requirement:
\[
{}^{x}\mathfrak{C}^i[\phi_j(x)] = \delta_{ij}, \quad \text{for } i,j = 1,2,\dots,k,
\]
where $\delta_{ij}$ denotes the Kronecker delta. This condition ensures that $\phi_j(x)$ remains neutral with respect to all other constraints while activating specifically at the $j^{\text{th}}$ constraint.
\end{itemize}

In this approach, original problem domain $[0,1]$ is transformed into standard interval $[-1,1]$ \cite{mortari2019high} to ensure compatibility with the definition of Chebyshev polynomials, which are employed as orthogonal basis functions for feature expansion. Since Chebyshev polynomials are inherently defined over the interval $[-1,1]$ and exhibit their orthogonality properties within this domain, such a transformation becomes necessary for their proper utilization.

To perform this transformation, a linear mapping is introduced that converts any point $X \in [0,1]$ to a corresponding point $X_i \in [-1,1]$. This mapping preserves the relative position of points while rescaling the interval appropriately. The transformation is given by 
\begin{equation*}
X_i = -1 + \frac{(1 - (-1))}{(1 - 0)} \, X,
\end{equation*}
which simplifies to
\begin{equation*}
X_i = -1 + 2X = 2X - 1.
\end{equation*}

Through this mapping, the left endpoint $X = 0$ is mapped to $X_i = -1$, and the right endpoint $X = 1$ is mapped to $X_i = 1$. Thus, the entire domain $[0,1]$ is linearly scaled onto $[-1,1]$, enabling the direct application of Chebyshev polynomials within the DFL-TFC framework.

Let us construct a CE that exactly satisfies all BCs. The CE of S-S boundary conditions and C-S boundary conditions are written respectively in mathematical form 

\begin{equation}\label{eq:constrained_expression}
   \boxed{\widehat{\overline{W}_p}(X) = h(X) + (X-1) h(0) - X h(1) + \frac{2X - 3X^2 + X^3}{6} h''(0) + \frac{X - X^3}{6} h''(1)  .}
 \end{equation}

\begin{equation}
\boxed{
\begin{aligned}  
 \widehat{\overline{W_p}}(X) = h(X) + \frac{-2 + 3X^2 - X^3}{2}h(0)+ \frac{-2X + 3X^2 - X^3}{2}h'(0) \\
 + \frac{X^3 - 3X^2}{2}h(1) + \frac{X^2 - X^3}{4}h''(1).
 \end{aligned}
 }
\end{equation}

  Detailed formulations of the CE for S-S and C-S boundary conditions are provided in the Appendix [\ref{appendix}].
Free function $h(X)$ in the CE is modeled using functional link neural network (FLNN), which maps input features $X$ into higher-dimensional space via Chebyshev polynomials of order 15, as this 15 order polynomial give better result as compare to 13 and 14 order (see the details in loss analysis Section \ref{loss analysis}), connects learnable weights to each basis function, takes the summation of all these terms, and passes the result through a linear activation function. Substituting $\widehat{\overline{W}_p}(X)$ into the \eqautoref{coleq21} yields the residual.

\begin{figure}[H]
    \centering
    \includegraphics[width=0.993\linewidth]{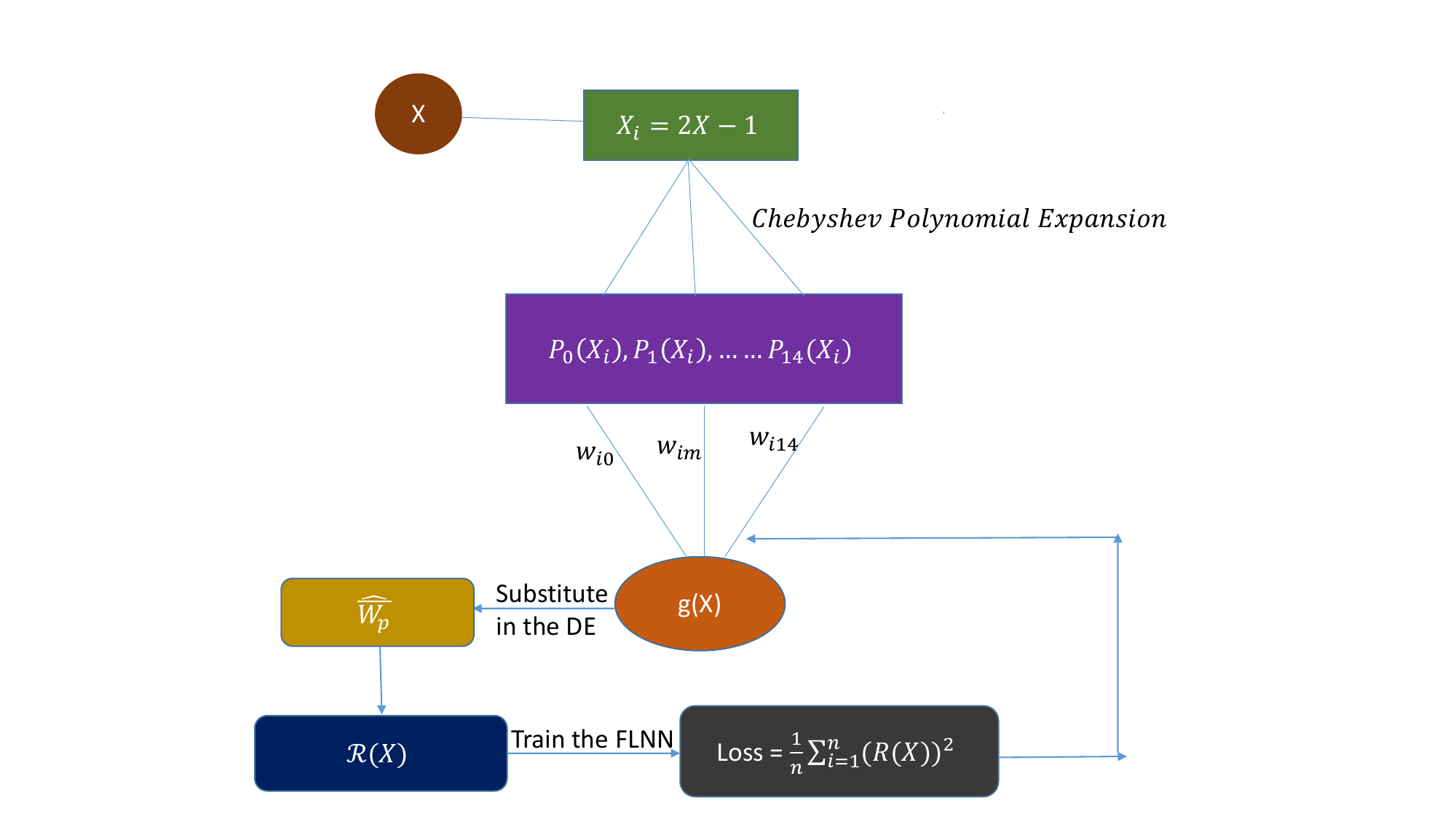}
    \caption{Schematic diagram of Domain mapped FL-TFC}
    \label{DFLTFC}
\end{figure}

Residual of the DE is given as:
\begin{equation*}
  \mathcal{R}(X) =   \frac{d^2}{dX^2} \left[ \widetilde{[EI]}_{eq} \frac{d^2 \widehat{\overline W_p}(X)}{dX^2} \right] -  \widetilde{ q_0} e^{\gamma X}- \widetilde{K_p} \frac{d^2 \widehat{\overline W_p}(X)}{dX^2}.
\end{equation*} 
For training purpose $100$ numbers of points are taken in the domain [0,1]. Gradient can be determined using automatic differentiation (AD)\cite{paszke2017automatic}. Loss is minimized via backpropagation using optimizers such as L-BFGS, executed in 10 sequential steps, where each step allows up to 50 iterations to achieve stable convergence and enhanced accuracy until the desired solution accuracy is achieved. Details of the above description are given in Fig. \ref{DFLTFC}.

Solving the above governing \eqautoref{coleq21} by the DFL-TFC method with both S-S and C-S boundary conditions, it will give the static deflection $\widehat{\overline W_p}(X)$. Then, for simplification, it is multiplied by 100, and it is represented as:
$$\widetilde{W_p}(X) = 100 \times \widehat{\overline W_p}(X), \  \ X \in [0,1].$$

\section{Results and Discussion}\label{sec. results}
In this section, we examine the validation of the DFL-TFC method and obtain new results, together with a loss analysis for both DFL-TFC and PINN.

\subsection{Validation of the DFL-TFC method} \label{sec. validation}
The proposed model, which describes the bending behavior of a tapered perforated beam while considering the foundation's effect with an exponential external load, was verified through a validation process in specific benchmark cases. In our problem, the Galerkin method (details given in \autoref{galerkin_sec}) is also used to obtain a reference solution for validation. This Galerkin framework serves as a reliable benchmark for evaluating the predictive capability of the DFL-TFC method, owing to its well-established accuracy in beam static bending studies. Since the analysis is carried out in a non-dimensional framework, all results are independent of specific material properties, and therefore, no particular beam material is prescribed. 

Since existing results for the bending deflection of tapered perforated beams are not available, the validation of the present model is performed in two stages. First, the maximum deflection of a simply supported (S-S) solid beam resting on an elastic foundation is compared with the results reported in the existing literature by Kafkas \cite{kafkas2025size} and Chen et al. \cite{chen2004mixed}. Table \ref{colT1} presents validation of the DFL-TFC results with existing literature for a solid beam by considering different values of foundation parameters $\widetilde{K_p}$ for maximum bending deflection. The accuracy of results is further supported by the Galerkin method.

\begin{figure}[H]
  \centering
  \begin{subfigure}[b]{0.49\linewidth}
    \includegraphics[width=\linewidth]{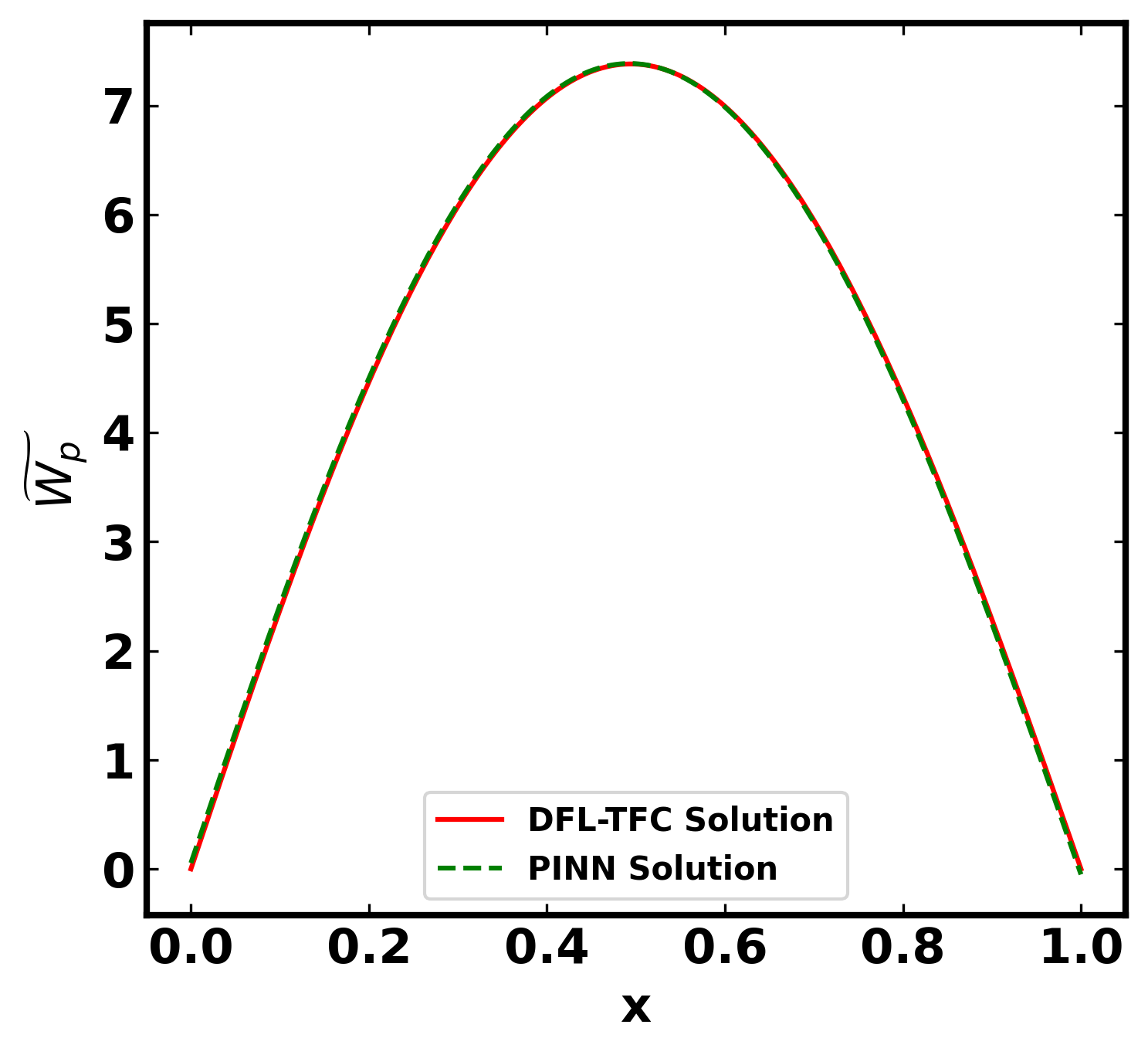}
     \caption{ For S-S boundary conditions}
     \label{val soln 1}
  \end{subfigure}
  \begin{subfigure}[b]{0.49\linewidth}
    \includegraphics[width=\linewidth]{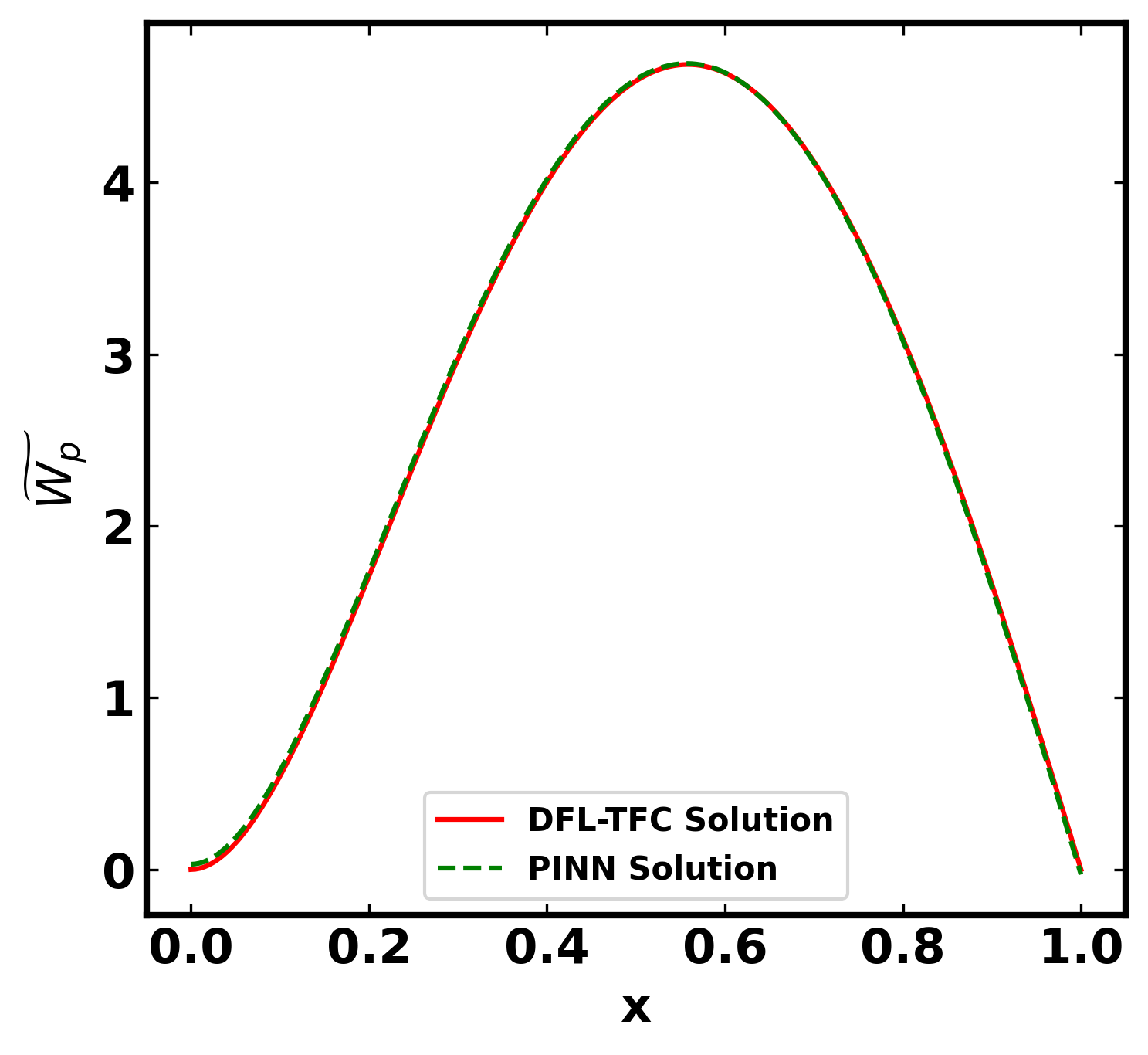}
    \caption{For C-S boundary conditions}
    \label{val soln 2}
  \end{subfigure}
  \caption{DFL-TFC and PINN solution for $\alpha=0.8, \gamma=5, \phi=0.5, \psi=0.5, N=3, \widetilde{K_p} =10, \widetilde{q_0}=1$}
  \label{validate solution}
\end{figure}

A close agreement between the present results and the existing results can be observed from Tables \ref{colT1}, \ref{colT2a}, \ref{colT3a} and Fig. \ref{validate solution}. This consistency indicates the accuracy of the proposed method and supports its reliability in predicting the bending response.

\begin{table}[H]
\centering
\caption{  Comparison of maximum bending deflection of a solid beam with uniform load $\widetilde{q_0}=1$  corresponding to $\alpha=1, \gamma=0, \phi=0, \psi=0, N=0 $}
\renewcommand{\arraystretch}{1.9}
\setlength{\tabcolsep}{9pt} 
\begin{tabular}{|c|c|c|c|c|c|}
\hline
Critical Deflection & \textbf{$\widetilde{K_p}$} & \textbf{DFL-TFC} & \textbf{ Galerkin} & \textbf{Existing} \cite{kafkas2025size} & \textbf{Existing} \cite{chen2004mixed}\\
 \hline
\multirow{3}{*}{$\widetilde{W}_p(0.5)$}
& $0$ & 1.3021 & 1.3021 & 1.3021 & 1.3023 \\
% \cline{2-6}
& $10$ & 0.6448 & 0.6448 & 0.6447 & 0.6448 \\
% \cline{2-6}
& $25$ &  0.3661 &  0.3661 & 0.3661 & 0.3661 \\
% \cline{2-6}
% & KP4 &  &  &  &  \\
\hline
\end{tabular}
\label{colT1}
\end{table}

Tables \ref{colT2a} and \ref{colT3a} provide a detailed comparison of the bending deflection of a tapered perforated beam, as predicted by the DFL-TFC, PINN, and Galerkin methods, for both S-S and C-S boundary conditions.

\begin{table}[H]
\centering
\caption{Comparison of bending deflection of a tapered perforated beam for both S-S and C-S boundary conditions corresponding to $\alpha = 0.3, N = 4, \gamma = 1, \phi = 0, \psi = 0, \widetilde{K_p} = 10, \text{and} \  \widetilde{q_0} = 10$}
\renewcommand{\arraystretch}{1.9}
 \setlength{\tabcolsep}{9pt} 
\begin{tabular}{|c|ccc|ccc|}
\hline
\multirow{2}{*}{$\widetilde{W}_p(X)$} & \multicolumn{3}{c|}{S-S} & \multicolumn{3}{c|}{C-S} \\
\cline{2-7}
 & DFL-TFC & PINN & Galerkin & DFL-TFC & PINN & Galerkin \\
\hline
$\widetilde{W}_p(0.1)$ & 4.4281 & 4.4280 & 4.4281 & 1.0552 & 1.1040 & 1.0552 \\
\hline
$\widetilde{W}_p(0.5)$ & 14.6425 & 14.6425 & 14.6425 & 9.5996 & 9.6123 & 9.5996 \\
\hline
$\widetilde{W}_p(0.9)$ & 5.0226 & 5.0227 & 5.0226 & 3.8835 & 3.8443 & 3.8835 \\
\hline

\end{tabular}
\label{colT2a}
\end{table}

\begin{table}[H]
\centering
\caption{Comparison of bending deflection of a tapered perforated beam for both S-S and C-S boundary conditions corresponding to $\alpha = 0.5, N = 2, \gamma = 0, \phi = 0.5, \psi = 0.5, \widetilde{K_p} = 10, \text{and} \  \widetilde{q_0} = 5$}
\renewcommand{\arraystretch}{1.9}
 \setlength{\tabcolsep}{9pt} 
\begin{tabular}{|c|ccc|ccc|}
\hline
\multirow{2}{*}{$\widetilde{W}_p(X)$} & \multicolumn{3}{c|}{S-S} & \multicolumn{3}{c|}{C-S} \\
\cline{2-7}
 & DFL-TFC & PINN & Galerkin & DFL-TFC & PINN & Galerkin \\
\hline
$\widetilde{W}_p(0.1)$ & 0.7980 & 0.8005 & 0.7980 & 0.1787 & 0.2090 & 0.1787 \\
\hline
$\widetilde{W}_p(0.5)$ & 2.1003 & 2.1001 & 2.1003 & 1.1714 & 1.1816 & 1.1714 \\
\hline
$\widetilde{W}_p(0.9)$ & 0.5610 & 0.5574 & 0.5610 & 0.3531 &  0.3313 & 0.3531 \\
\hline

\end{tabular}
\label{colT3a}
\end{table}

Next, the bending deflection of a tapered perforated beam subjected to an exponential load is examined by comparing the present results obtained from our proposed DFL-TFC method with those computed using the traditional PINN approach. This comparison is carried out to assess the accuracy and reliability of our proposed technique. 
In Fig. \ref{validate solution}, the solutions obtained using DFL-TFC and PINN are presented for the parameters
$\alpha = 0.8,\ \gamma = 5,\ \phi = 0.5,\ \psi = 0.5,\ N = 3,\ \widetilde{K_p} = 10,\ \text{and} \  \widetilde{q_0} = 1$ considering both S-S and C-S boundary conditions.
 For the S-S case, the DFL-TFC method achieves a mean square loss of \(5.3257\times 10^{-11}\), while the PINN method reaches \(4.7997\times 10^{-6}\). For the C-S case, the DFL-TFC mean square loss reduces to \(3.7958\times 10^{-11}\), and where as PINN mean square loss is \(3.0377\times 10^{-5}\).

\subsection{New Outcomes}
The bending deflection of a tapered perforated beam is obtained using the DFL-TFC method, which demonstrates faster convergence, higher accuracy, and greater computational efficiency compared with the traditional physics-informed neural network (PINN) approach. Table \ref{colT2} presents bending deflection behavior of the tapered perforated beam under simply supported-simply supported (S-S) and clamped-simply supported (C-S) boundary conditions for different values of filling ratio ($\alpha$) and number of rows of holes ($N$). For the present analysis, the following parameter values are considered: $\gamma = 1$, $\phi = 0$, $\psi = 0$, $\widetilde{K_p} = 10$, and $\widetilde{q_0} = 10$.

From Table \ref{colT2}, it can be observed that bending deflection decreases with an increase in filling ratio ($\alpha$). This behavior occurs because a higher filling ratio indicates that more material is retained in the beam, which enhances its stiffness and consequently reduces the deflection. On the other hand, as the number of rows of holes ($N$) increases, deflection tends to increase. This is mainly due to additional perforations that reduce the effective stiffness of the beam, leading to a larger bending response.

\begin{table}[H]
\centering
\caption{Bending deflection of a tapered perforated beam for C-S boundary conditions with varying $\alpha$ and $N$ corresponding to $\gamma = 1, \phi = 0, \psi = 0, \widetilde{K_p} = 10, \text{and} \  \widetilde{q_0} = 10$} 

\renewcommand{\arraystretch}{1.9}
 \setlength{\tabcolsep}{13pt} 
\begin{tabular}{|c|c|ccc|ccc|}
\hline
\multirow{2}{*}{$\alpha$} & \multirow{2}{*}{$N$}
& \multicolumn{3}{c|}{S-S}
& \multicolumn{3}{c|}{C-S} \\
\cline{3-8}
 &  & $\widetilde{W}_p(0.1)$ & $\widetilde{W}_p(0.5)$ & $\widetilde{W}_p(0.9)$
 & $\widetilde{W}_p(0.1)$ & $\widetilde{W}_p(0.5)$ & $\widetilde{W}_p(0.9)$ \\
 \hline

\multirow{2}{*}{0.1}
& 1 & 4.4966 &  14.8698 &  5.1111 & 1.0881 & 9.8393 & 3.9801 \\
 % \cline{2-8}
& 4 &  5.5348  & 18.2472 &  6.5745 & 1.8192 & 14.1231 & 5.7176 \\
% \hline

\multirow{2}{*}{0.3}
& 1 & 3.5462 & 11.7104 & 3.9349 & 0.7134 & 6.8760 & 2.7865 \\
% \cline{2-8}
& 4 & 4.4281 & 14.6425 & 5.0226 & 1.0552 & 9.5996 & 3.8835 \\
% \hline

\multirow{2}{*}{0.5}
& 1 & 3.3416 & 11.0303 & 3.6930 & 0.6500 & 6.3260 & 2.5647\\
% \cline{2-8}
& 4 & 3.8069 & 12.5772 & 4.2480 & 0.8015 & 7.6175 & 3.0854 \\
% \hline

\multirow{2}{*}{0.7}
& 1 & 3.3002 & 10.8928 & 3.6444 & 0.6377 & 6.2180 & 2.5211 \\
% \cline{2-8}
& 4 & 3.4607 & 11.4262 & 3.8335 & 0.6863 & 6.6429 & 2.6926 \\
% \hline

\multirow{2}{*}{0.9}
& 1 & 3.2957 &  10.8778 & 3.6392 & 0.6364 & 6.2063 & 2.5164 \\
% \cline{2-8}
& 4 & 3.3121 & 10.9323 & 3.6584 & 0.6413 & 6.2489 & 2.5336 \\
\hline

\end{tabular}
\label{colT2}
\end{table}

The variation of bending deflection in the tapered perforated beam under S-S and C-S boundary conditions, corresponding to different values of the tapered parameters $\phi$ and $\psi$, is summarized in Table \ref{ColT3}. The analysis is conducted with the following set of parameters: $\gamma = 0, \alpha =0.5, N=2, \widetilde{K_p}=10, \text{and} \  \widetilde{q_0} =5$

From Table \ref{ColT3}, it can be observed that bending deflection decreases as the tapered parameters $\phi$ and $\psi$ increase for both boundary conditions. This behavior occurs because larger values of the tapering parameters modify beam mass by effectively increasing beam stiffness. As a result, the beam becomes more resistant to bending under the applied load. Consequently, the overall deflection is reduced for both S-S and C-S boundary conditions.

\begin{table}[H]
\centering
\caption{Bending deflection of a tapered perforated beam for S-S and C-S boundary conditions with varying  $\phi$ and $\psi$ corresponding to 
$\gamma = 0, \alpha =0.5, N=2, \widetilde{K_p}=10, \text{and} \  \widetilde{q_0} =5$}

\renewcommand{\arraystretch}{1.9}
\setlength{\tabcolsep}{13pt}

\begin{tabular}{|c|c|ccc|ccc|}
\hline
$\phi$ & $\psi$
& \multicolumn{3}{c|}{S-S}
& \multicolumn{3}{c|}{C-S} \\
\cline{3-8}

& & $\widetilde{W}_p(0.1)$ & $\widetilde{W}_p(0.5)$ & $\widetilde{W}_p(0.9)$
& $\widetilde{W}_p(0.1)$ & $\widetilde{W}_p(0.5)$ & $\widetilde{W}_p(0.9)$ \\
\hline

\multirow{3}{*}{0.1}
& 0.1 & 1.0450 & 3.1588 & 0.9587 & 0.2135 & 1.7805 & 0.6363\\
& 0.5 & 0.9652 & 2.7368 & 0.7662 & 0.2015 & 1.5275 & 0.4918\\
& 0.9 & 0.8959 & 2.3948 & 0.6326 & 0.1909 & 1.3251 & 0.3954\\

% \hline
\multirow{3}{*}{0.5}
& 0.1 &  0.8612 &  2.4030 & 0.6811 & 0.1891 & 1.3494 & 0.4392\\
& 0.5 & 0.7980 &  2.1003 & 0.5610 & 0.1787 & 1.1714 & 0.3531\\
& 0.9 & 0.7438 & 1.8561 & 0.4742 & 0.1697 & 1.0289 & 0.2927\\

% \hline
\multirow{3}{*}{0.9}
& 0.1 & 0.7108 & 1.8439 & 0.4979 & 0.1671 & 1.0389 & 0.3164\\
& 0.5 & 0.6626 & 1.6310 & 0.4209 & 0.1586 & 0.9144 & 0.2626\\
& 0.9 & 0.6211 & 1.4574 & 0.3631 & 0.1511 & 0.8133 & 0.2230\\

\hline
\end{tabular}
\label{ColT3}
\end{table}

The bending deflection of a tapered perforated beam is evaluated for S-S and C-S boundary conditions. The effects of the exponential load parameter $\gamma$ and the foundation parameter $\widetilde{K}_p$ are presented in Table \ref{ColT4}. In the present study, the parameters are taken as $\alpha = 0.8$, $N = 3$, $\phi = \psi = 0.5$, and $\widetilde{q_0} = 1$.

As seen in Table~\ref{ColT4}, the deflection of the tapered perforated beam increases with an increase in the exponential load parameter $\gamma$. This is because higher values of $\gamma$ correspond to a stronger external loading on the beam, resulting in larger bending deformation. Conversely, increasing the foundation parameter $\widetilde{K}_p$ provides greater support stiffness to the beam, which restricts its deformation and leads to smaller deflection under both boundary conditions.

\begin{table}[H]
\centering
\caption{Bending deflection of a tapered perforated beam for S-S and C-S boundary conditions with varying  $\gamma$ and $\widetilde{K_p}$ corresponding to 
$ \alpha =0.8, N=3, \phi = \psi =0.5, \text{and} \  \widetilde{q_0} =1$}

\renewcommand{\arraystretch}{1.9}
\setlength{\tabcolsep}{13pt}

\begin{tabular}{|c|c|ccc|ccc|}
\hline
$\gamma$ & $\widetilde{K_p}$
& \multicolumn{3}{c|}{S-S}
& \multicolumn{3}{c|}{C-S} \\
\cline{3-8}

& & $\widetilde{W}_p(0.1)$ & $\widetilde{W}_p(0.5)$ & $\widetilde{W}_p(0.9)$
& $\widetilde{W}_p(0.1)$ & $\widetilde{W}_p(0.5)$ & $\widetilde{W}_p(0.9)$ \\
\hline

\multirow{3}{*}{1}
& 1 & 0.3235 & 0.8699 & 0.2342 & 0.0592 & 0.4255 & 0.1322\\
& 5 & 0.2741 & 0.7442 & 0.2023 & 0.0552 & 0.3927 & 0.1225\\
& 10 & 0.2298 & 0.6307 & 0.1734 & 0.0509 & 0.3582 & 0.1122\\

% \hline
\multirow{3}{*}{3}
& 1 &  0.9575 &  2.7292 & 0.7730 & 0.1789 & 0.1789 & 0.4723\\
& 5 & 0.8069 &  2.3391 & 0.6729 & 0.1658 & 1.3098 & 0.4390\\
& 10 & 0.6723 & 1.9865 & 0.5817 & 0.1520 & 1.1939 & 0.4038\\

% \hline
\multirow{3}{*}{5}
& 1 & 3.4081 & 10.1141 & 2.9981 & 0.6443 & 5.4663 & 1.9310\\
& 5 & 2.8609 & 8.6782 & 2.6264 & 0.5949 & 5.0401 & 1.7998\\
& 10 & 2.3729 & 7.3792 & 2.2867 & 0.5431 & 4.5918 & 1.6609\\

\hline
\end{tabular}
\label{ColT4}
\end{table}

\subsubsection{Effect of filling ratio}
\begin{figure}[H]
  \centering
  \begin{subfigure}[b]{0.49\linewidth}
    \includegraphics[width=\linewidth]{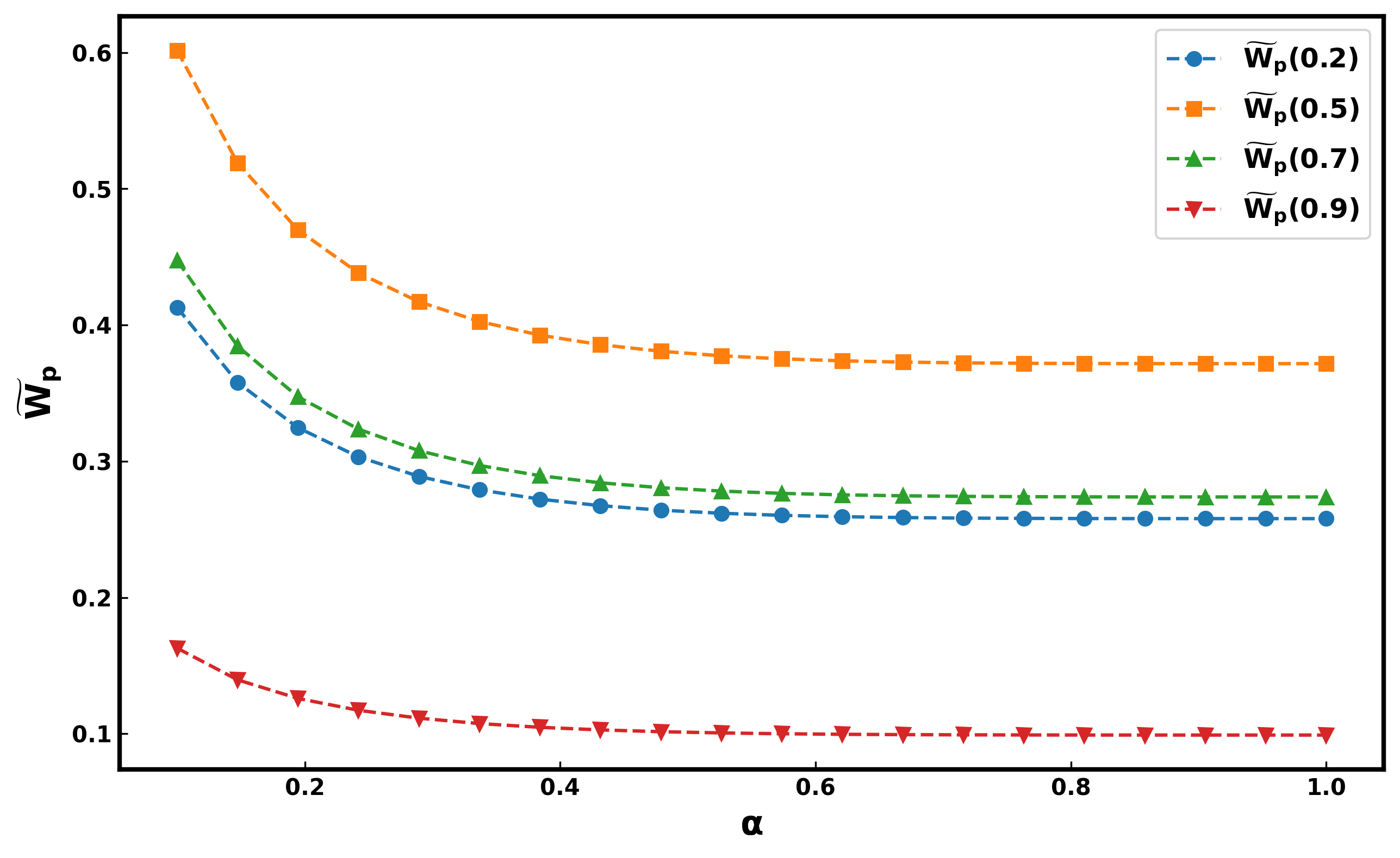 }
     \caption{S-S tapered perforated beam}
     \label{fig. alpha1}
  \end{subfigure}
  \begin{subfigure}[b]{0.49\linewidth}
    \includegraphics[width=\linewidth]{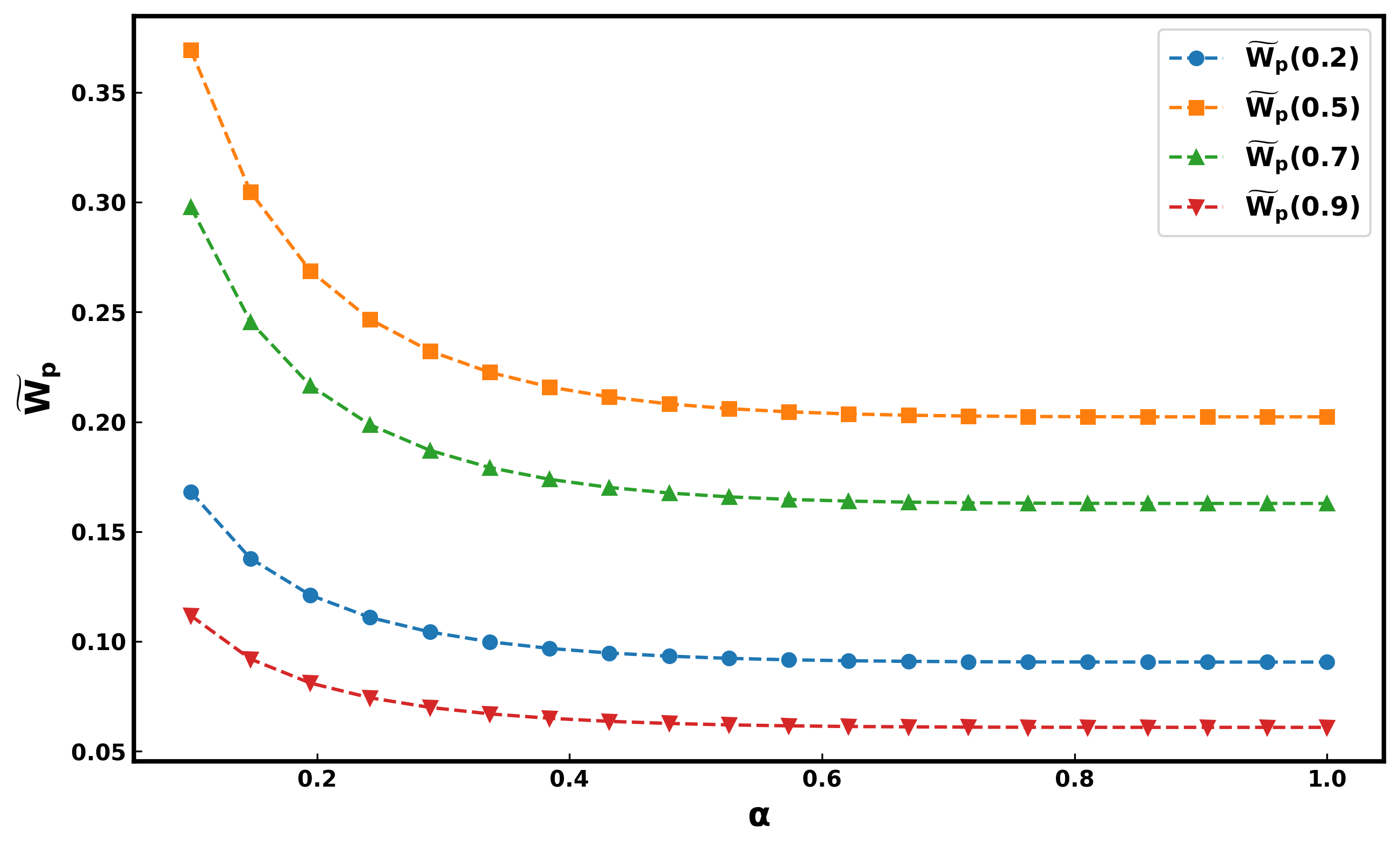}
    \caption{C-S tapered perforated beam}
    \label{fig. alpha2}
  \end{subfigure}
  \caption{Variation of bending deflection with respect to filling ratio $\alpha$ for a tapered perforated beam subjected to an exponential load}
  \label{colfig4}
\end{figure}

The variation of bending deflection in the tapered perforated beam with respect to filling ratio $(\alpha)$ under exponential loading and in the presence of an elastic foundation with parameter $(\widetilde{K_p})$ is shown in Fig. \ref{colfig4}. The bending behavior is evaluated for both S-S and C-S boundary conditions. Fig. \ref{colfig4} presents bending deflection as the filling ratio varies from 0 to 1 for specified $\widetilde{W_p}$  at $X=0.2,0.5, 0.7, \text{and}\  0.9$. The analysis is conducted with the following set of parameters: $N=1$, $\phi = \psi = 0.5$, $\gamma = 0$, $\widetilde {q_0} = 1$, and $\widetilde {K_p} =10$.

From Fig. \ref{colfig4}, it can be seen that deflection gradually decreases as $\alpha$ increases from 0 to 1 for both S-S and C-S boundary conditions. This trend occurs because an increase in the filling ratio reduces effective perforation size, thereby increasing the solid portion of the beam. As a consequence, the stiffness of the tapered perforated beam becomes higher, which restricts the bending deformation.

\subsubsection{Effect of exponential load parameter}

\begin{figure}[H]
  \centering
  \begin{subfigure}[b]{0.49\linewidth}
    \includegraphics[width=\linewidth]{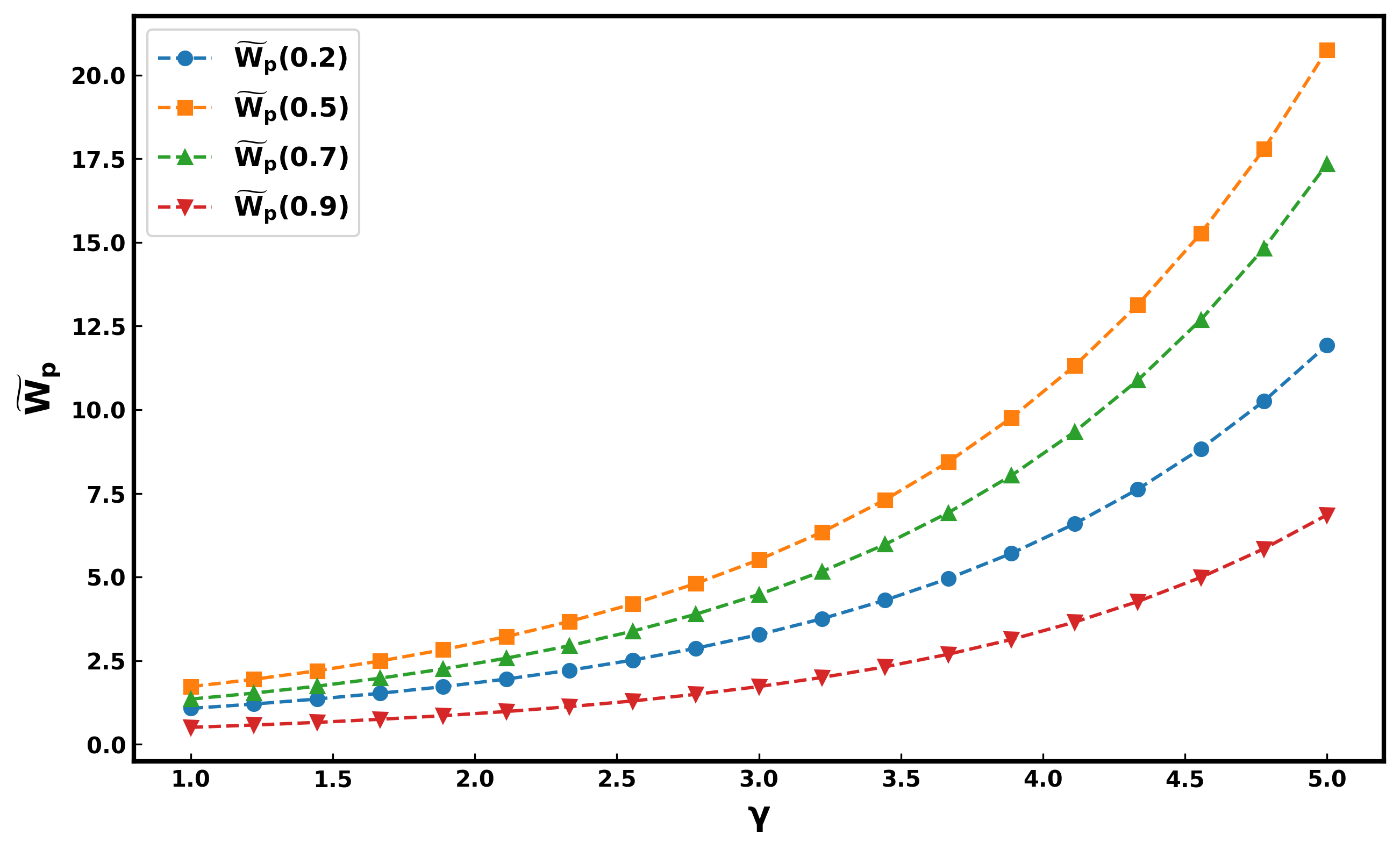}
     \caption{S-S tapered perforated beam }
     \label{fig. gamma1}
  \end{subfigure}
  \begin{subfigure}[b]{0.49\linewidth}
    \includegraphics[width=\linewidth]{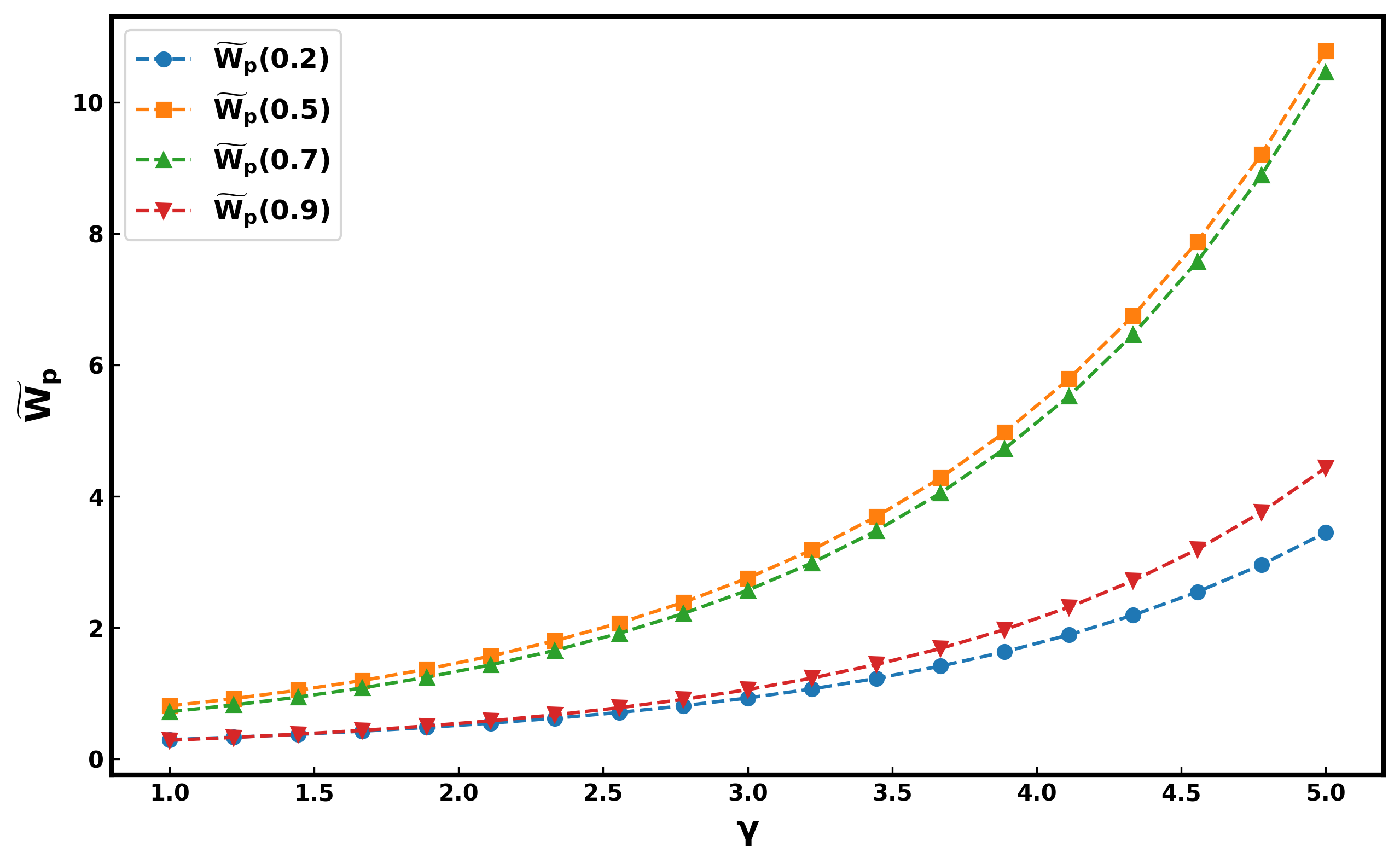}
    \caption{C-S tapered perforated beam}
    \label{fig. gamma2}
  \end{subfigure}
  \caption{Variation of bending deflection with respect to exponential load parameter $\gamma$ for a tapered perforated beam resting on foundation}
   \label{colfig5}
\end{figure}

Fig. \ref{colfig5} presents the bending deflection response of a tapered perforated beam for varying values of the exponential load parameter $\gamma$ in presence of an elastic foundation. The parameter $\gamma$ is varied from 1 to 5 in the present analysis. The corresponding bending deflections are evaluated at the locations $X = 0.2$, $0.5$, $0.7$, and $0.9$. The results are presented for both S-S and C-S boundary conditions. The parameter values used in this analysis are as follows: $\alpha=0.5$, $N=3$, $\phi = \psi = 0.2$, $\widetilde {q_0} = 1$, and $\widetilde {K_p} =1$.

As the exponential load parameter $\gamma$ increases, the bending deflection of the tapered perforated beam also increases, as clearly shown in Fig. \ref{colfig5}. Higher values of $\gamma$ enhance the load acting along the beam, which results in greater bending. Since the tapered beam contains perforations, the presence of holes reduces the overall stiffness of the beam, making it more sensitive to changes in the applied load. Consequently, larger deflections are observed for both S-S and C-S boundary conditions.

\subsubsection{Effect of tapering parameters}
\begin{figure}[H]
  \centering
  \begin{subfigure}[b]{0.49\linewidth}
    \includegraphics[width=\linewidth]{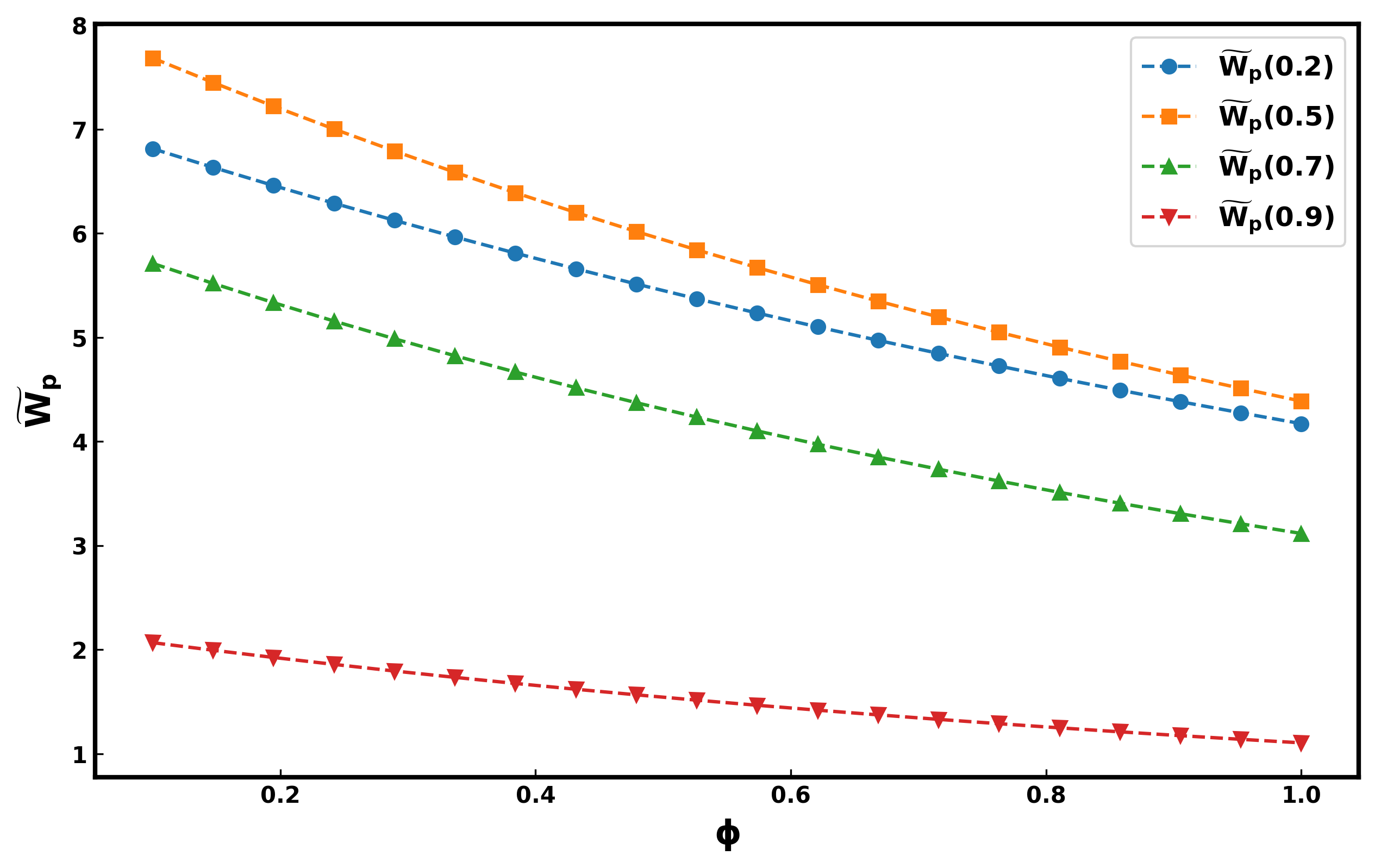 }
     \caption{S-S tapered perforated beam}
     \label{fig. phi1}
  \end{subfigure}
  \begin{subfigure}[b]{0.49\linewidth}
    \includegraphics[width=\linewidth]{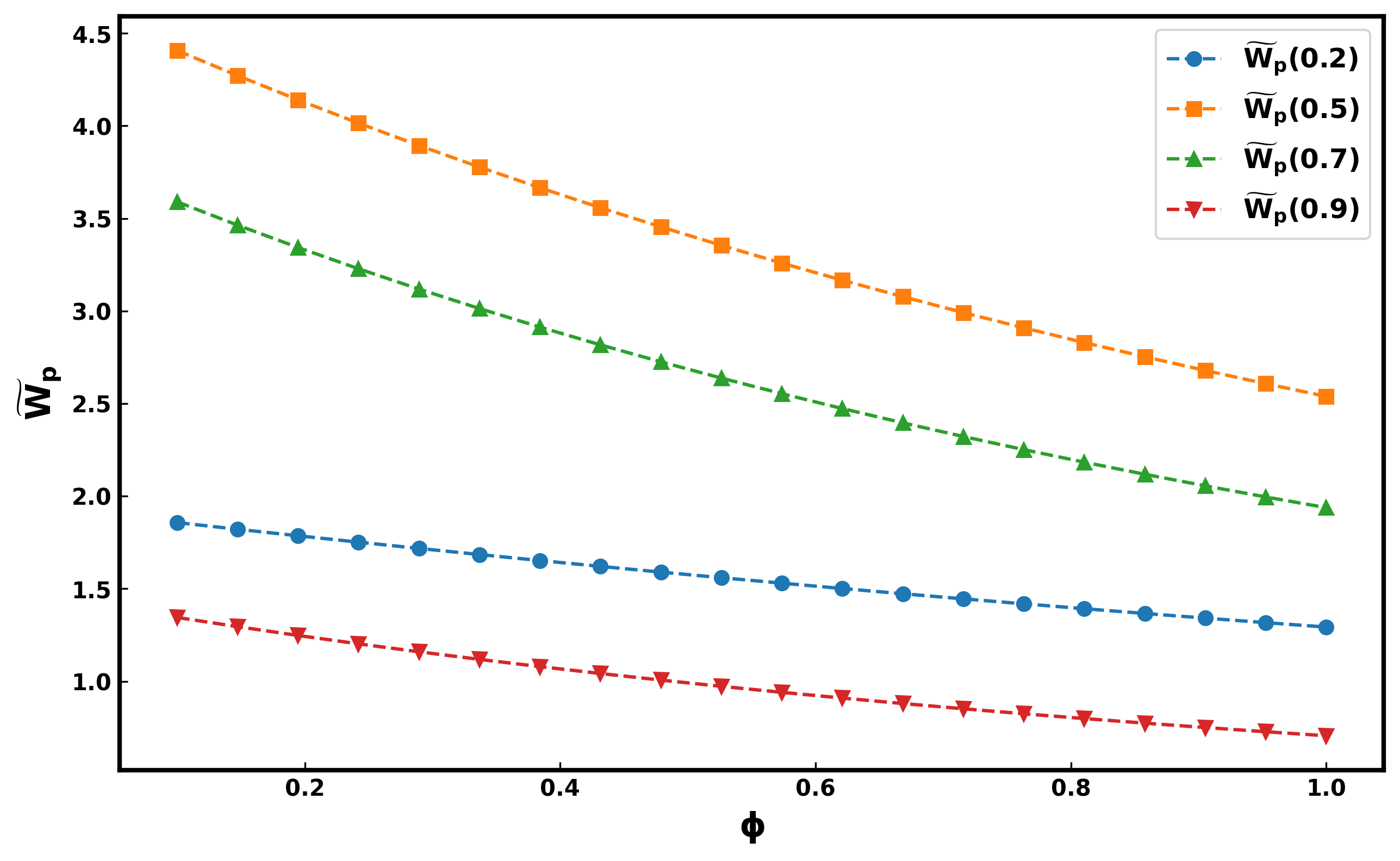}
    \caption{C-S tapered perforated beam}
    \label{fig. phi2}
  \end{subfigure}
  \caption{Effect of tapering parameter $\phi$ on the bending behaviour of a tapered perforated beam resting on a foundation}
  \label{colfig6}
\end{figure}

\begin{figure}[H]
  \centering
  \begin{subfigure}[b]{0.49\linewidth}
    \includegraphics[width=\linewidth]{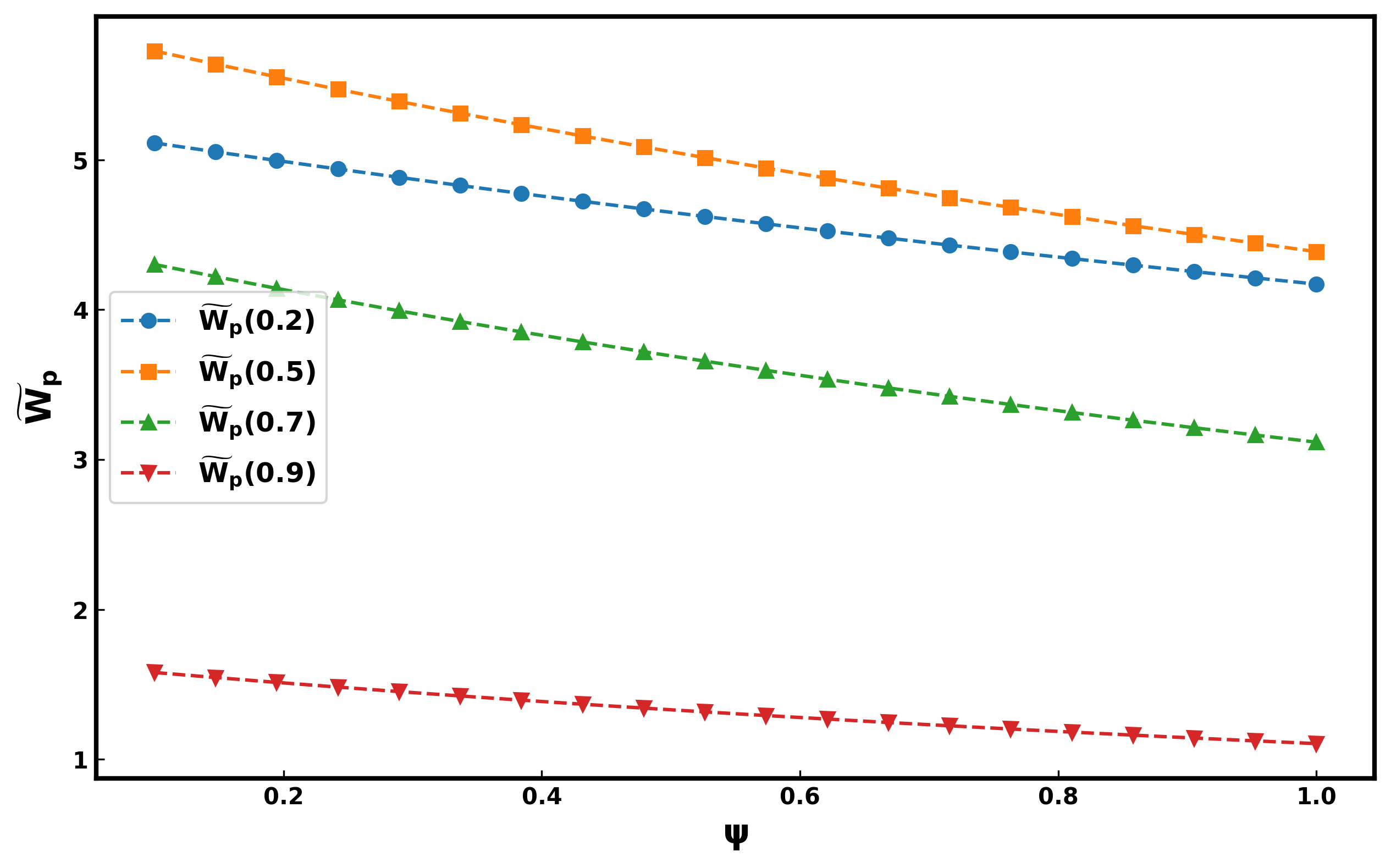 }
     \caption{S-S tapered perforated beam}
     \label{fig. psi1}
  \end{subfigure}
  \begin{subfigure}[b]{0.49\linewidth}
    \includegraphics[width=\linewidth]{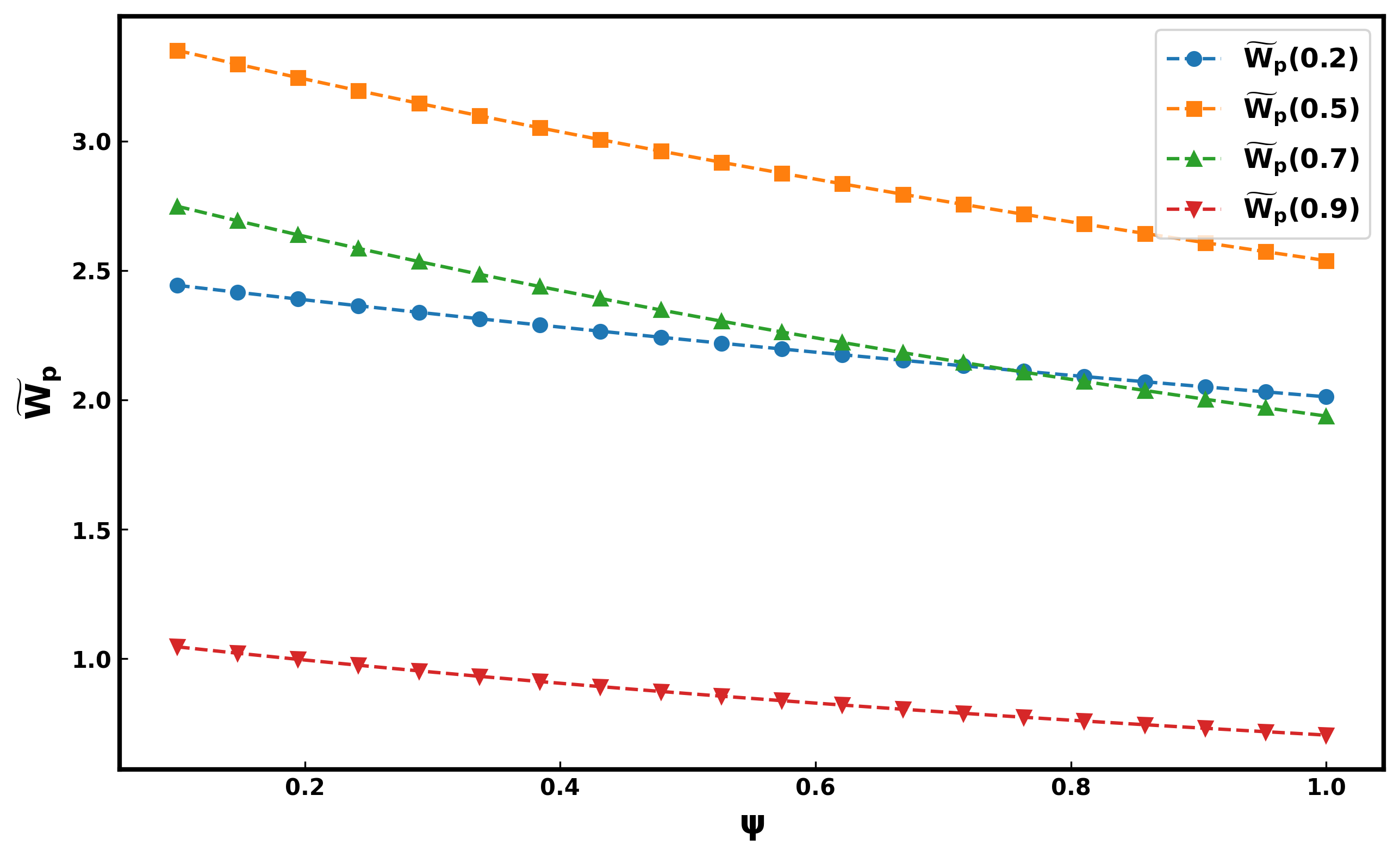}
    \caption{C-S tapered perforated beam}
    \label{fig. psi2}
  \end{subfigure}
  \caption{Effect of tapering parameter $\psi$ on the bending behaviour of a tapered perforated beam resting on a foundation}
  \label{colfig7}
\end{figure}

It is worth mentioning that the tapering parameters $\phi$ and $\psi$ play an important role in determining the bending deflection of a tapered perforated beam resting on an elastic foundation under exponential loading. Their effect on deflection behavior is illustrated in Figs. \ref{colfig6} and \ref{colfig7}. The results are presented for both S-S and C-S boundary conditions. These figures highlight how variations in tapering parameters modify deflection characteristics of the beam.

Two different sets of parameters are considered in this analysis:

For Fig. \ref{colfig6}: $\alpha=0.5$, $N=2$, $ \psi = 0.2$, $\widetilde {q_0} = 10$, and $\widetilde {K_p} =10$.

For Fig. \ref{colfig7}: $\alpha=0.5$, $N=2$, $ \phi = 0.2$, $\widetilde {q_0} = 10$, and $\widetilde {K_p} =10$.

A clear reduction in bending deflection can be observed in Fig. \ref{colfig6} as tapering parameters $\phi$ increase. This parameter modifies cross-sectional dimensions of the beam through linear variations along its length, effectively adding more material to the perforated beam. As a result, the overall stiffness of the tapered perforated beam becomes higher. The increased stiffness limits the bending deformation and leads to smaller deflection values. The same decreasing trend is evident for both S-S and C-S boundary conditions. A similar trend can also be observed in Fig. \ref{colfig7}, consistent with the results shown in Fig. \ref{colfig6}.

\subsubsection{Effect of foundation parameter}
\begin{figure}[H]
  \centering
  \begin{subfigure}[b]{0.49\linewidth}
    \includegraphics[width=\linewidth]{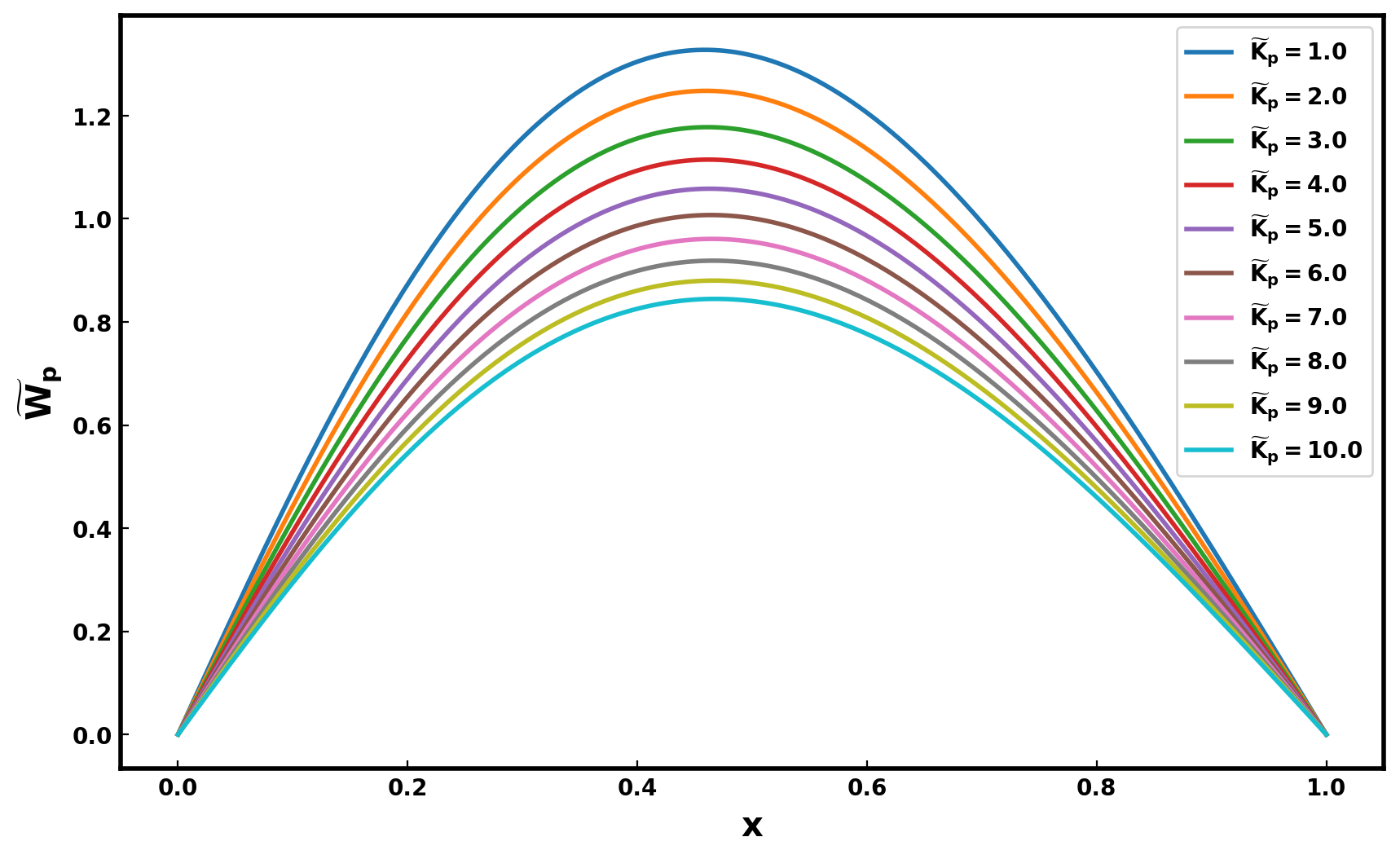 }
     \caption{S-S tapered perforated beam }
     \label{fig. w_pred for s-s bc}
  \end{subfigure}
  \begin{subfigure}[b]{0.49\linewidth}
    \includegraphics[width=\linewidth]{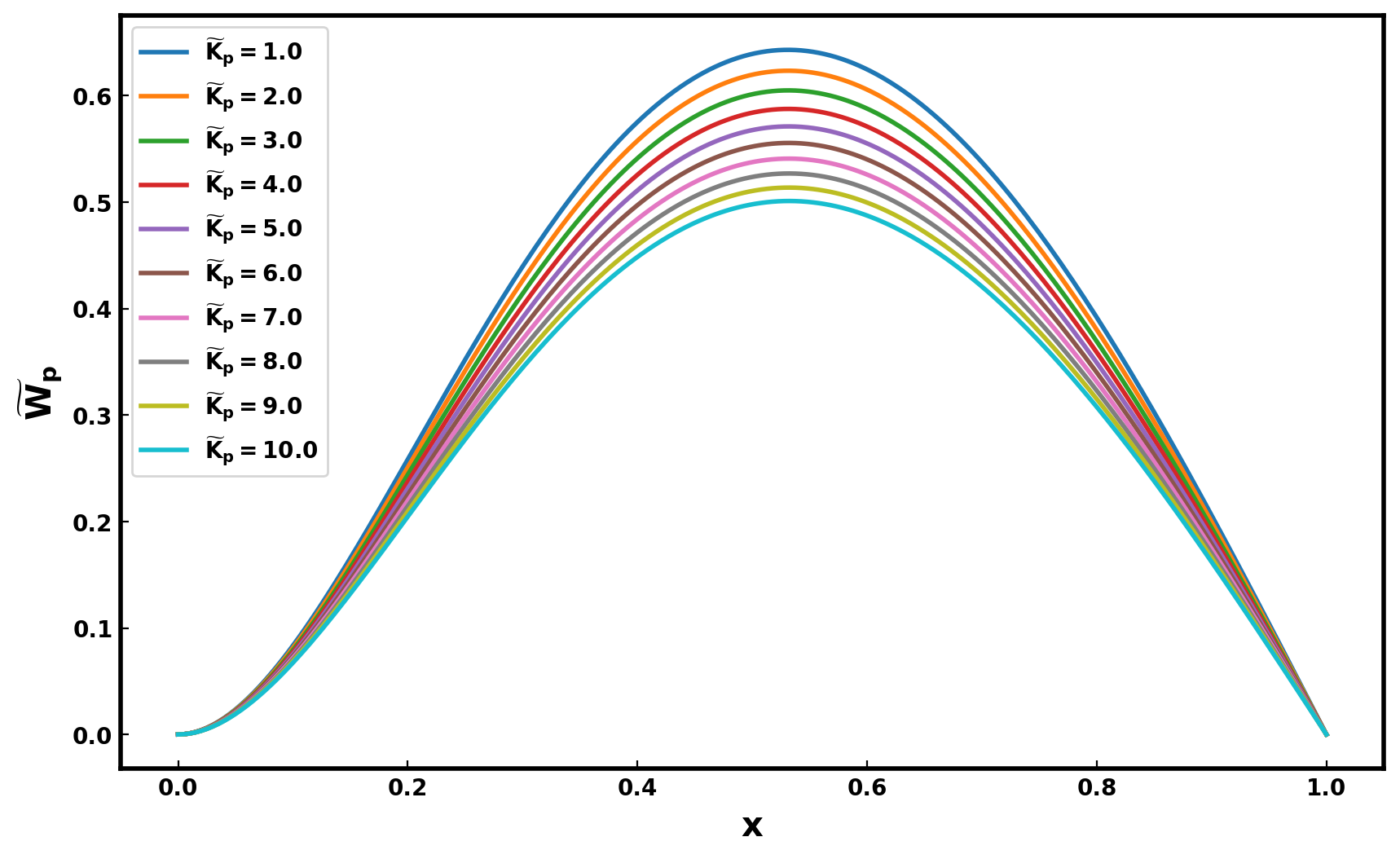}
    \caption{C-S tapered perforated beam}
    \label{fig. w_pred for c-s bc}
  \end{subfigure}
  \caption{Effect of elastic foundation parameter $\widetilde{K_p}$ on the bending behaviour of a tapered perforated beam with exponential load}
   \label{colfig8}
\end{figure}

It is important to mention that the foundation parameter $\widetilde{K_p}$ strongly influences the bending deflection of the tapered perforated beam. The effect of this  foundation parameter $\widetilde{K_p}$ is presented in Fig. \ref{colfig8}  for both of S-S and C-S boundary conditions. The range of this parameter is 1 to 10 for bending deflection. The following parameter values are taken for conducting the analysis: $\alpha=0.5$, $N=4$, $\phi = \psi = 0.4$, $\widetilde {q_0} = 1$, and $\gamma =1$.

A decreasing trend in bending deflection can be observed in Fig. \ref{colfig8}, as the foundation parameter $\widetilde{K_p}$ increases for both boundary conditions. A higher value of $\widetilde{K_p}$ represents a stiffer elastic foundation. This additional stiffness provides greater support to the beam and restricts its bending under the applied load. As a result, overall deflection of the tapered perforated beam becomes smaller.

\subsection{Loss Analysis of  DFL-TFC and PINN }\label{loss analysis}

In this section, we present a detailed discussion on loss behavior of domain-mapped FL-TFC and PINN approaches under both S-S and C-S boundary conditions. Training dataset used in all cases consists of 100 points distributed across the domain. For all results, we employ the LBFGS optimizer, configured with 50 optimization steps in PINN and 10 steps in DFL-TFC, where each step allows for a maximum of 50 internal iterations to ensure sufficient convergence.

To analyze the effect of network architecture in the PINN framework, we consider models with one, two, and three hidden layers, where each hidden layer contains five neurons with tanh activation function in the hidden layer, while output layer uses linear activation function. This configuration enables a systematic investigation of how network depth influences convergence behavior and predictive performance while keeping the number of neurons per layer fixed.
For DFL-TFC method, loss behavior is examined by varying order of Chebyshev polynomial expansion used in functional approximation. Specifically, polynomial orders 13, 14, and 15 are considered to assess the impact of increasing the approximation order on solution accuracy and convergence of the loss function. Results demonstrate how higher-order polynomial expansions enhance approximation capability and contribute to improved convergence characteristics. 

\begin{figure}[H]
  \centering
  \begin{subfigure}[b]{0.30\linewidth}
    \includegraphics[width=\linewidth]{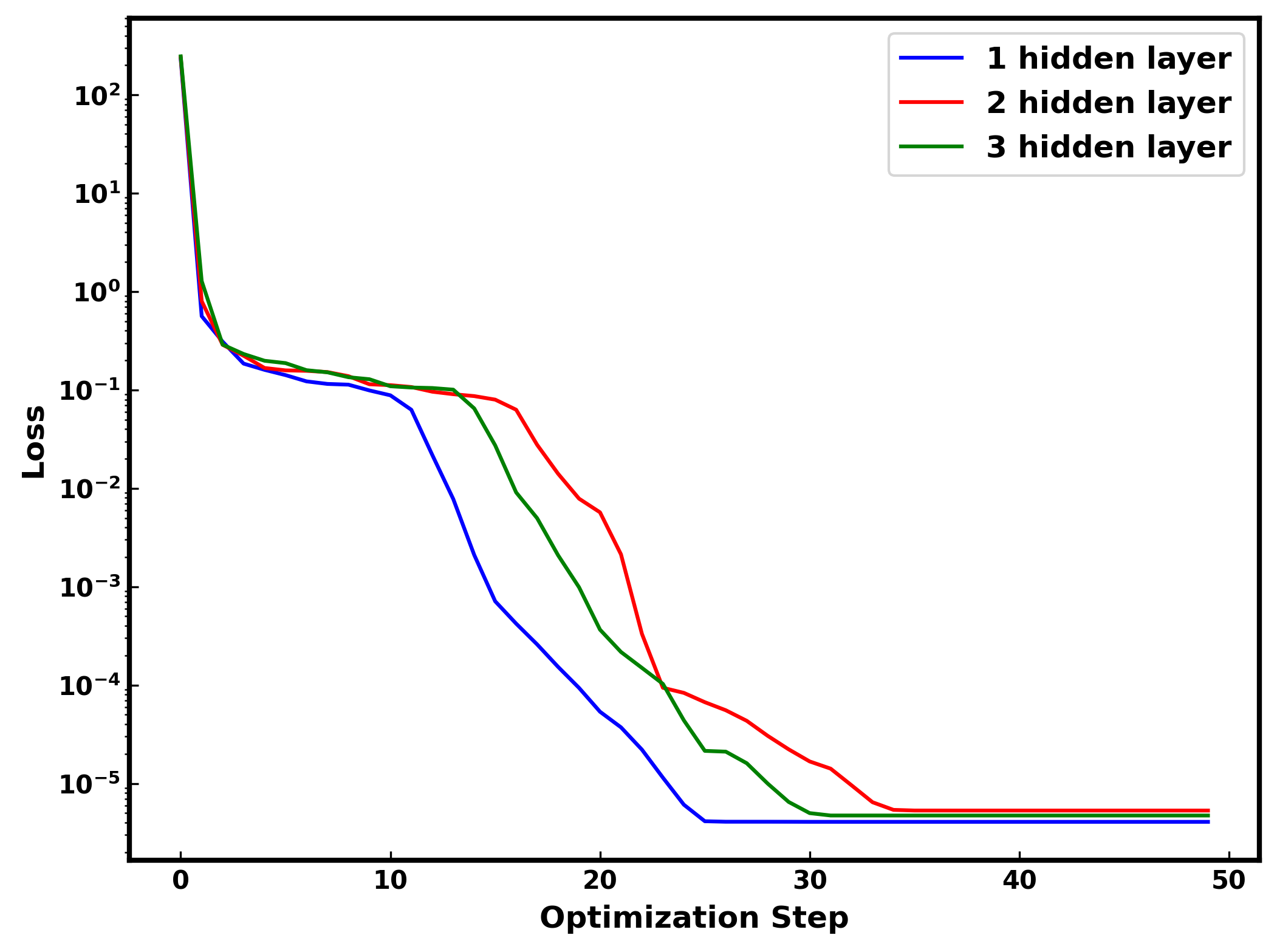}
     % \caption{$\alpha, \gamma, \phi, \psi, N, k_p, q_0 = 0.7, 3, 0.5, 0.5, 1, 10, 2$}
     \caption{Loss 1}
     \label{loss fig1}
  \end{subfigure}
  \begin{subfigure}[b]{0.30\linewidth}
    \includegraphics[width=\linewidth]{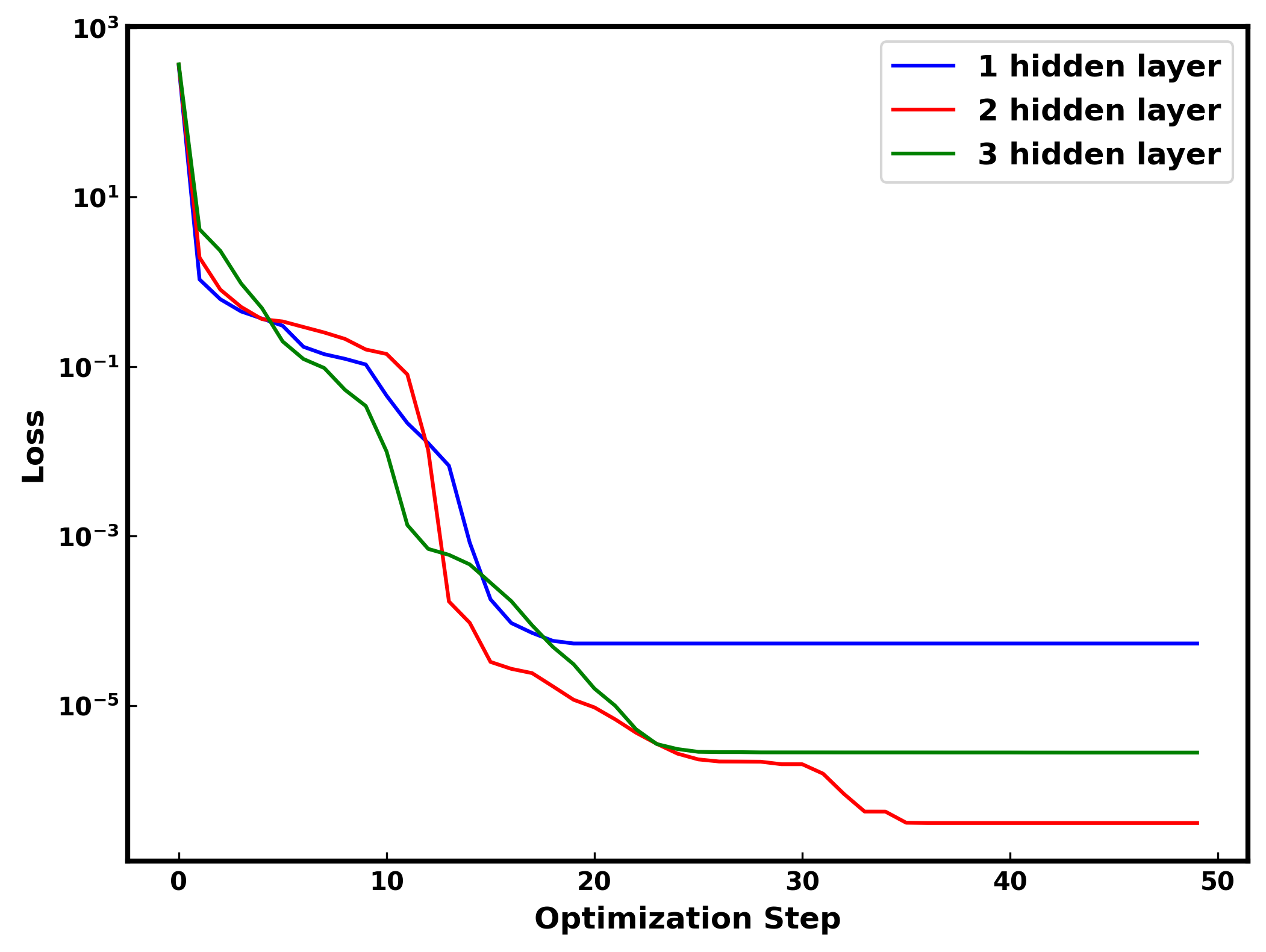}
    % \caption{$alpha, gamma, phi, psi, N_val, kp, q0 = 0.8, 4, 0.3, 0.3, 2, 5, 1$}
    \caption{Loss 2}
    \label{loss fig2}
  \end{subfigure}
   \begin{subfigure}[b]{0.30\linewidth}
    \includegraphics[width=\linewidth]{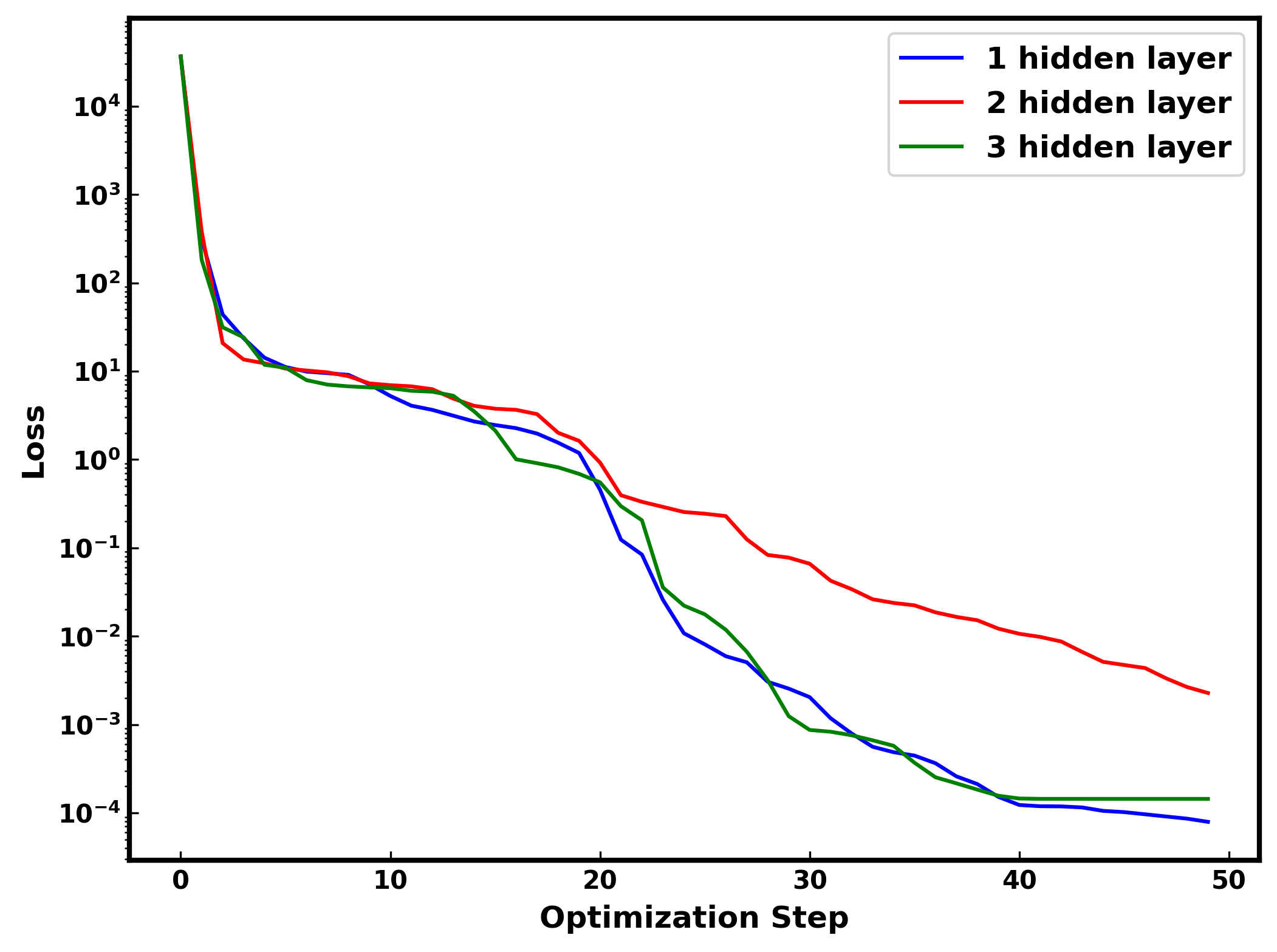}
    % \caption{$\alpha, \gamma, \phi, \psi, N, k_p, q_0 = 0.2, 5, 0.2, 0.2, 1, 8, 4$}
    \caption{Loss 3}
    \label{loss fig3}
  \end{subfigure}
   \begin{subfigure}[b]{0.30\linewidth}
    \includegraphics[width=\linewidth]{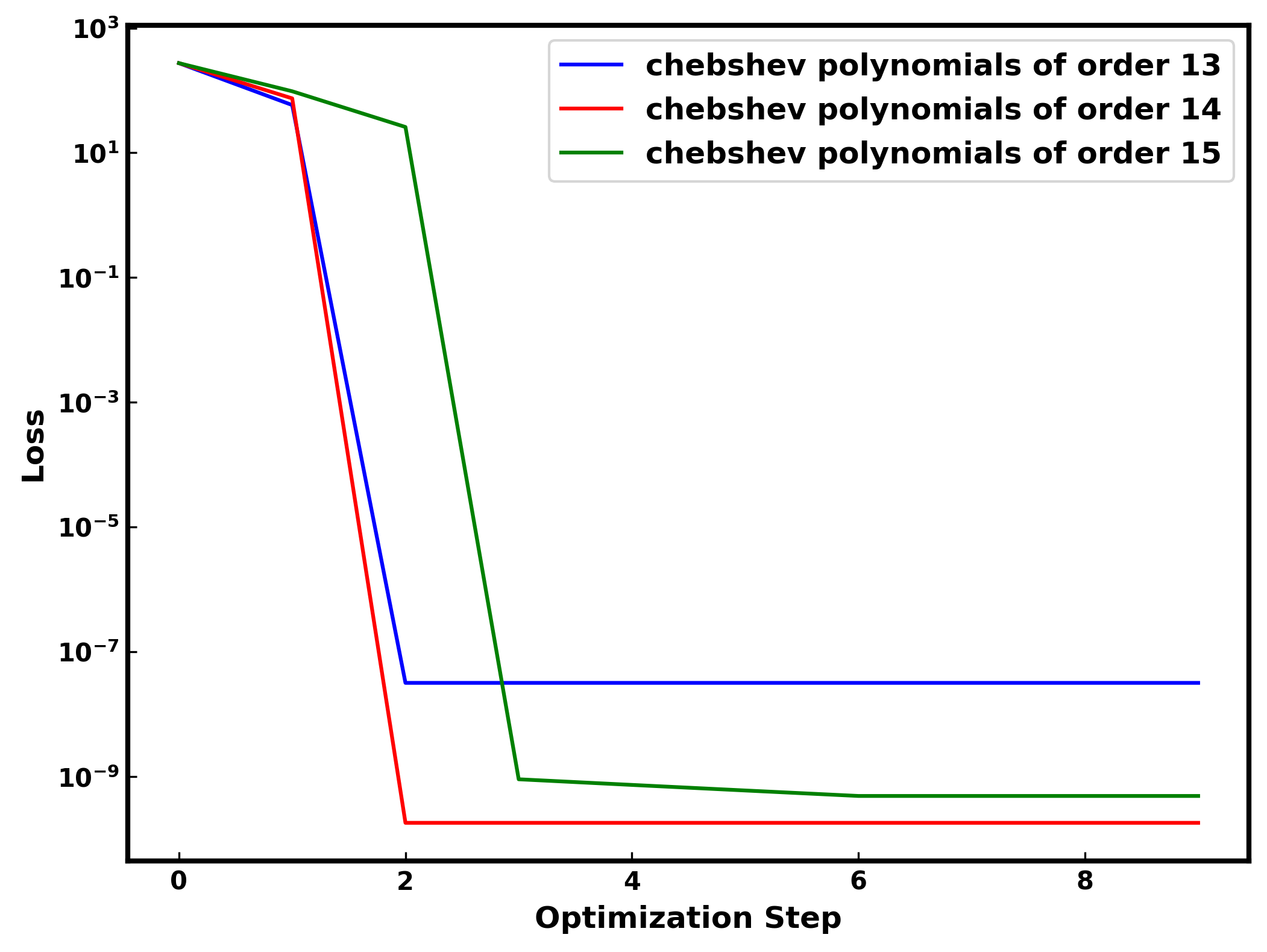}
    \caption{Loss 4}
    \label{loss fig4}
  \end{subfigure}
  \begin{subfigure}[b]{0.30\linewidth}
    \includegraphics[width=\linewidth]{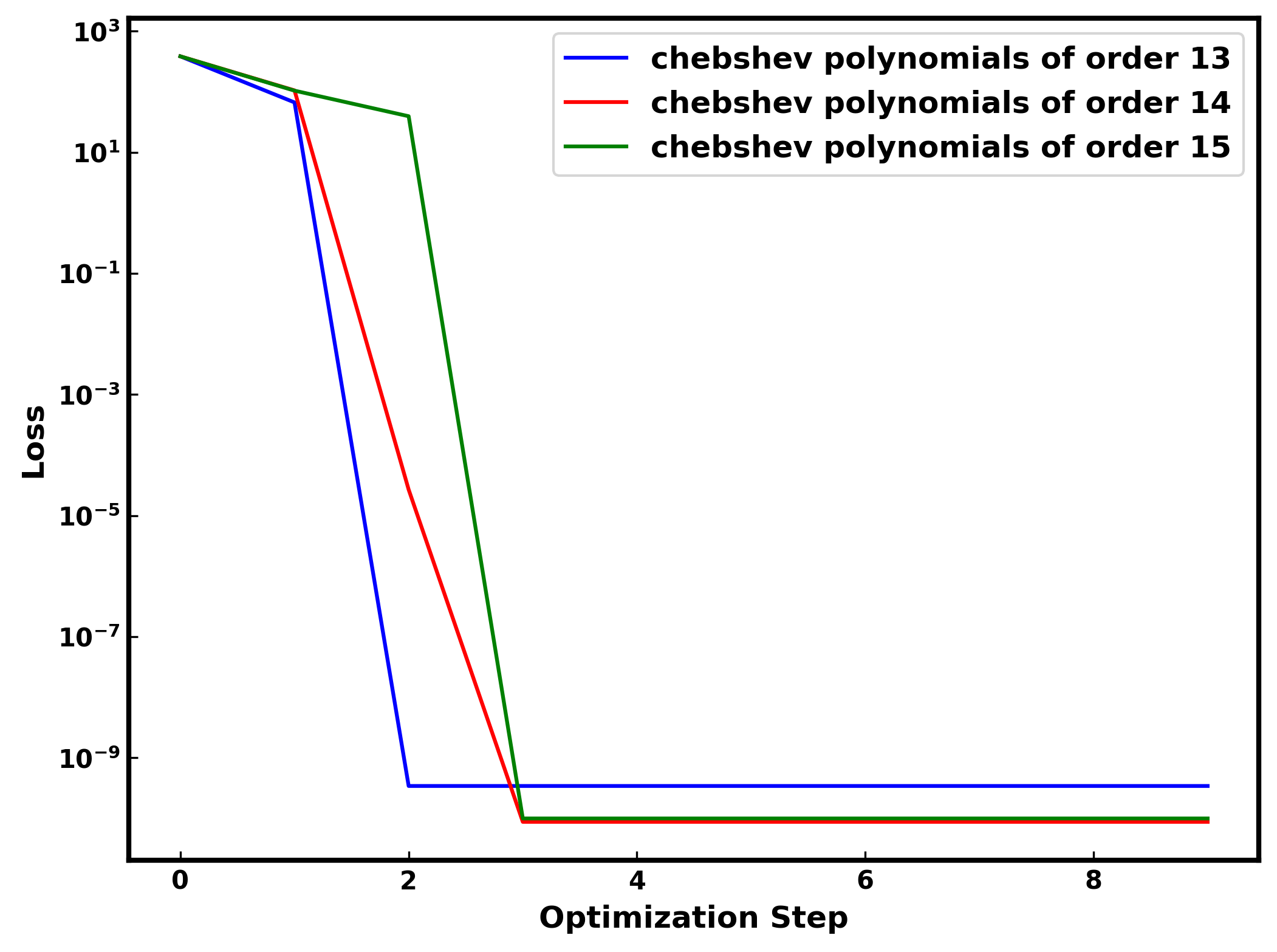}
    \caption{Loss 5}
    \label{loss fig5}
  \end{subfigure}
   \begin{subfigure}[b]{0.30\linewidth}
    \includegraphics[width=\linewidth]{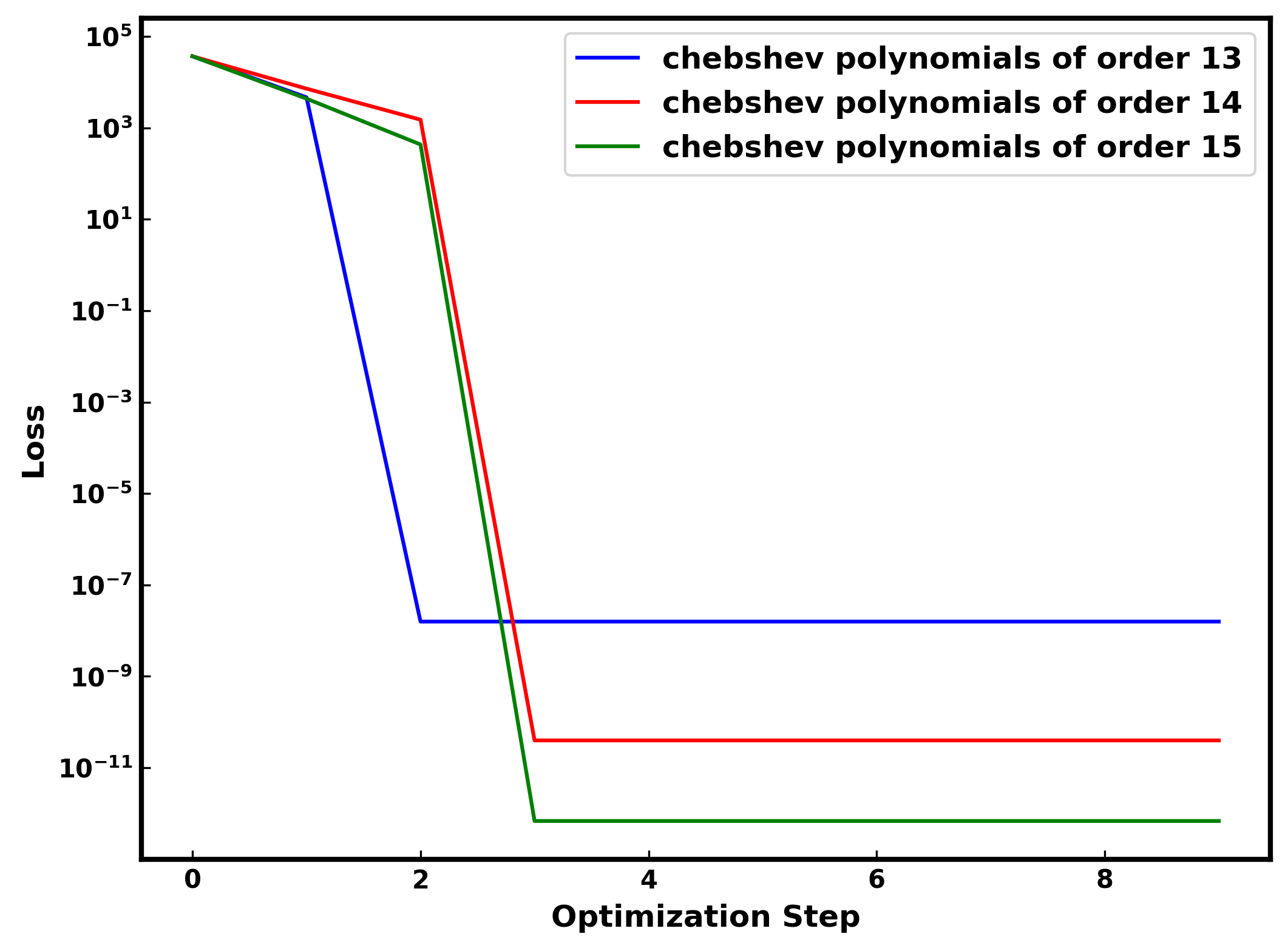}
    \caption{Loss 6}
    \label{loss fig6}
  \end{subfigure}
  \caption{Loss figures of PINN and DFL-TFC frameworks at different parameters for S-S boundary conditions }
   \label{loss figure 1}
\end{figure}

The details corresponding to Figs \ref{loss fig1} and \ref{loss fig4} are presented in Table \ref{comparison table 1}. Likewise, the details of Figs \ref{loss fig2} and \ref{loss fig5} are given in Table \ref{comparison table 2}, while those of Figs \ref{loss fig3} and \ref{loss fig6} are summarized in Table \ref{comparison table 3} for S-S boundary conditions.

\begin{figure}[H]
  \centering
  \begin{subfigure}[b]{0.30\linewidth}
    \includegraphics[width=\linewidth]{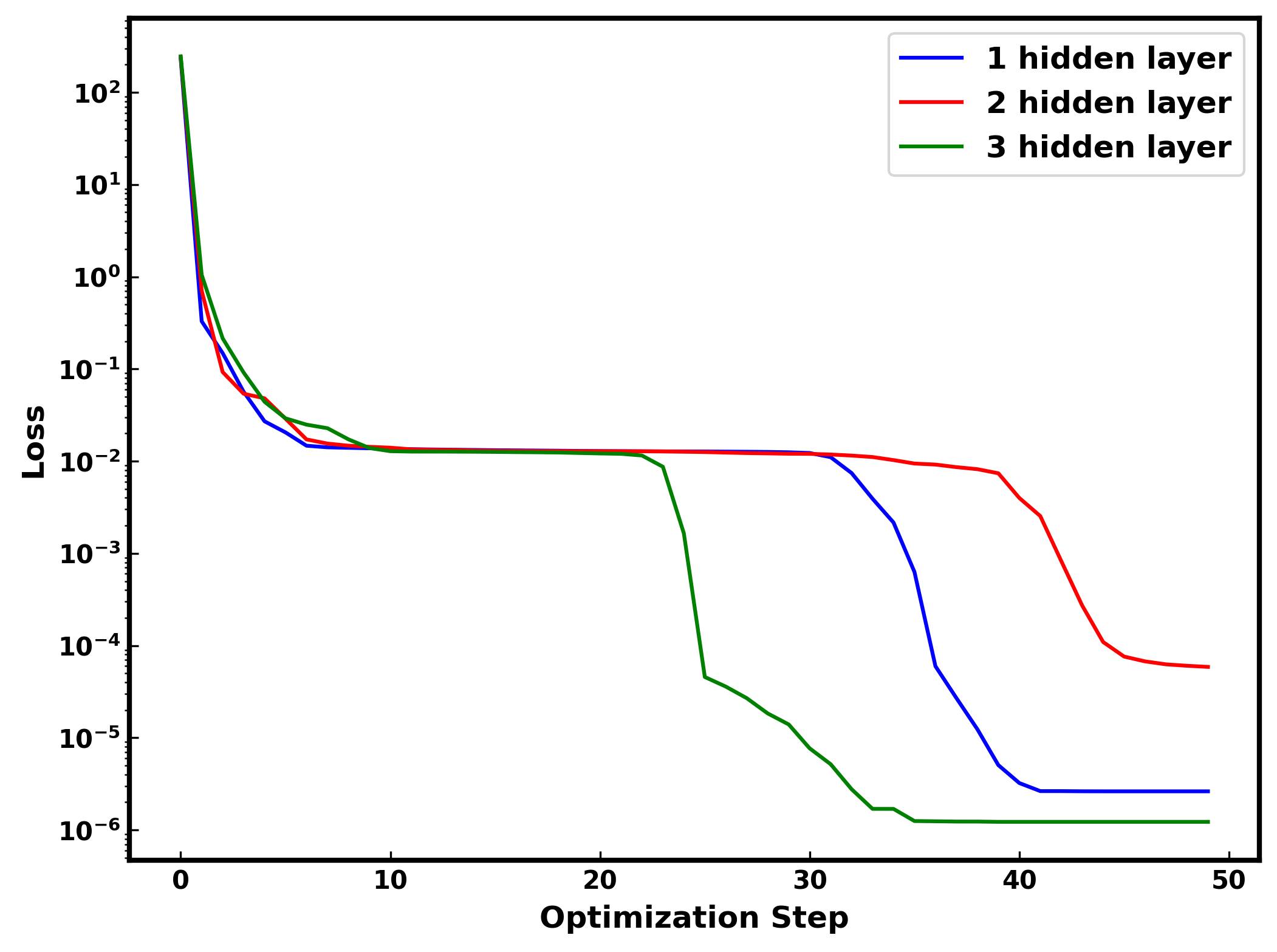}
     % \caption{$\alpha, \gamma, \phi, \psi, N, k_p, q_0 = 0.7, 3, 0.5, 0.5, 1, 10, 2$}
     \caption{Loss 7}
     \label{loss fig7}
  \end{subfigure}
  \begin{subfigure}[b]{0.30\linewidth}
    \includegraphics[width=\linewidth]{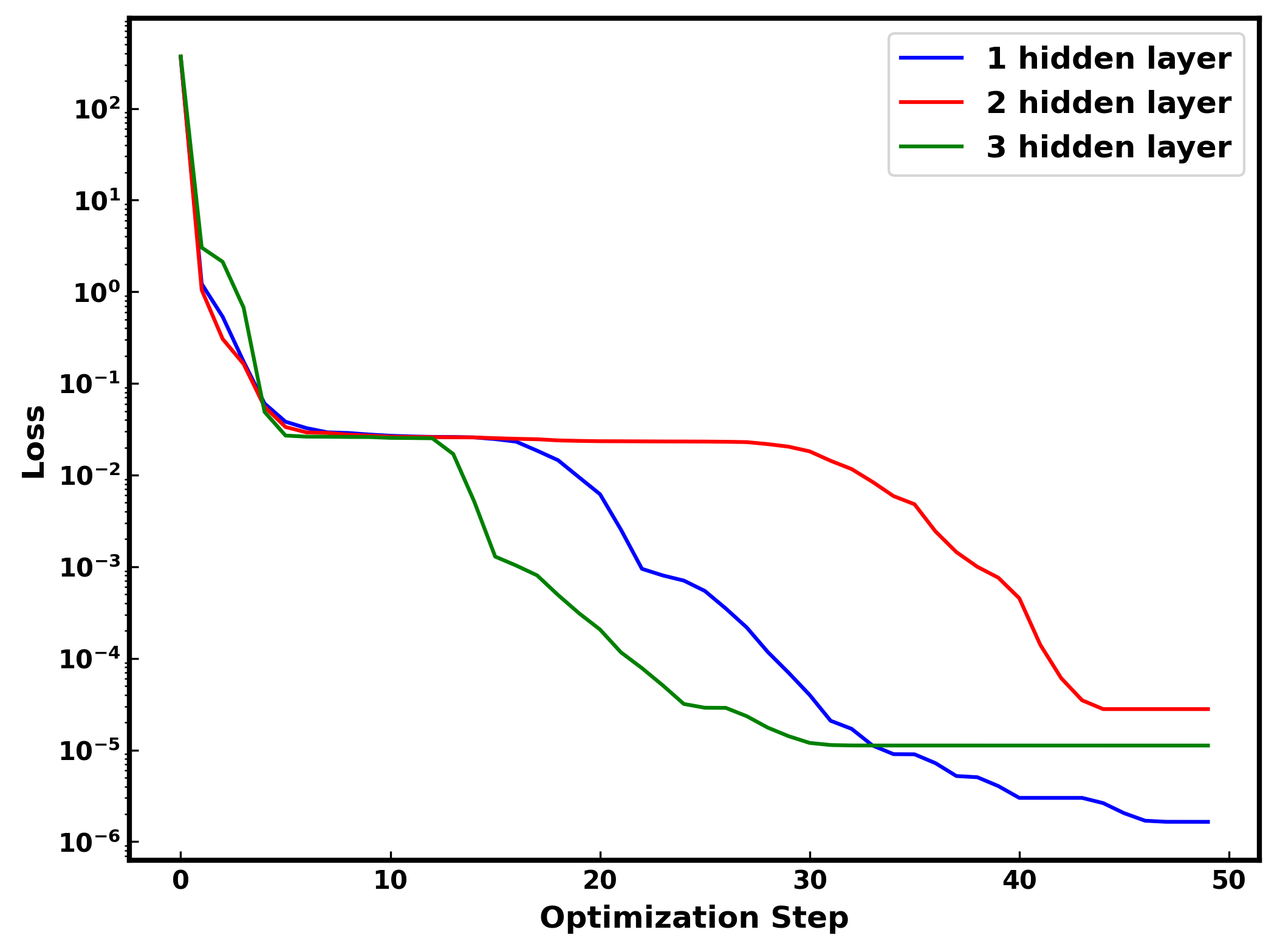}
    % \caption{$alpha, gamma, phi, psi, N_val, kp, q0 = 0.8, 4, 0.3, 0.3, 2, 5, 1$}
    \caption{Loss 8}
    \label{loss fig8}
  \end{subfigure}
   \begin{subfigure}[b]{0.30\linewidth}
    \includegraphics[width=\linewidth]{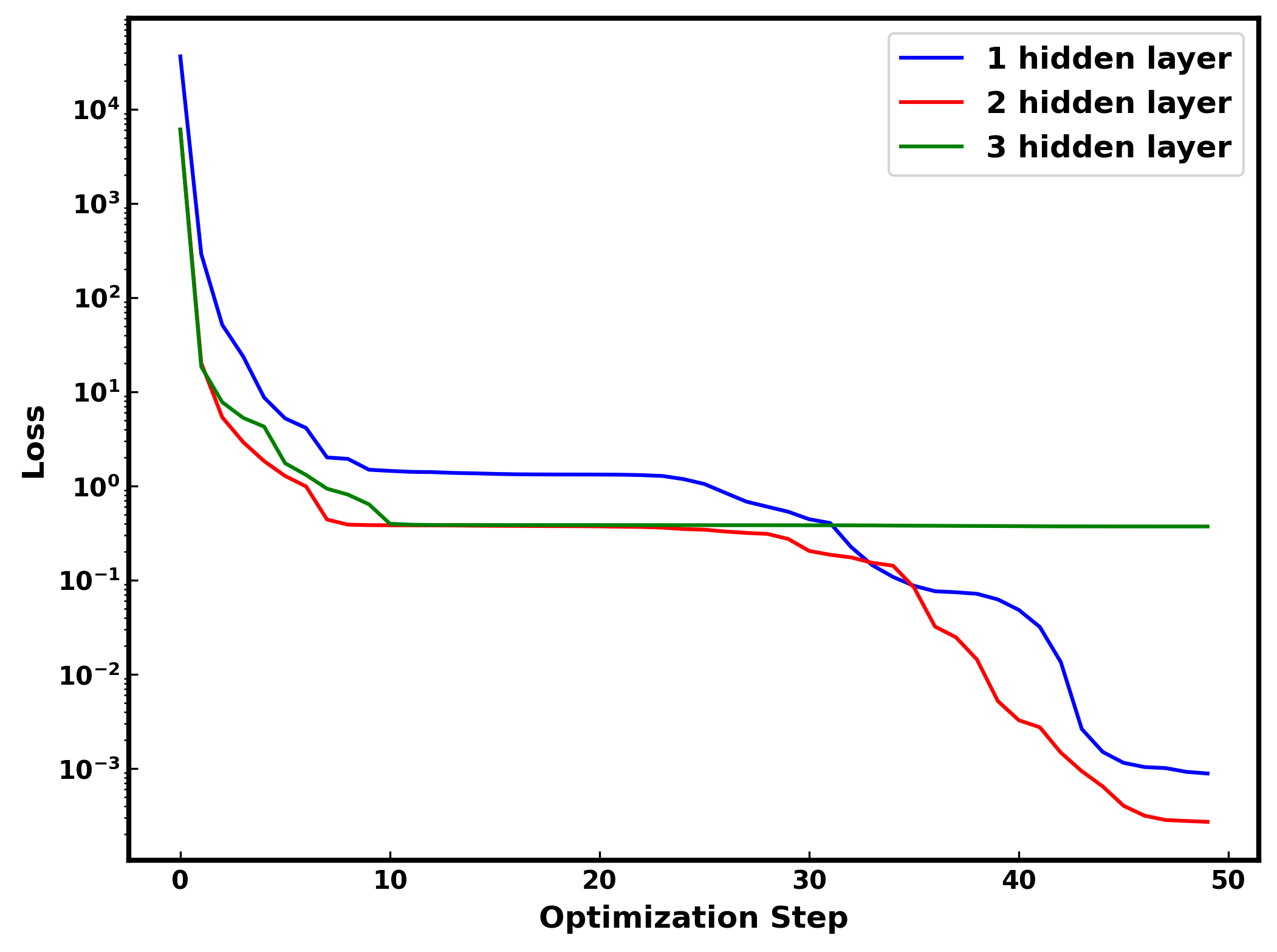}
    % \caption{$\alpha, \gamma, \phi, \psi, N, k_p, q_0 = 0.2, 5, 0.2, 0.2, 1, 8, 4$}
    \caption{Loss 9}
    \label{loss fig9}
  \end{subfigure}
   \begin{subfigure}[b]{0.30\linewidth}
    \includegraphics[width=\linewidth]{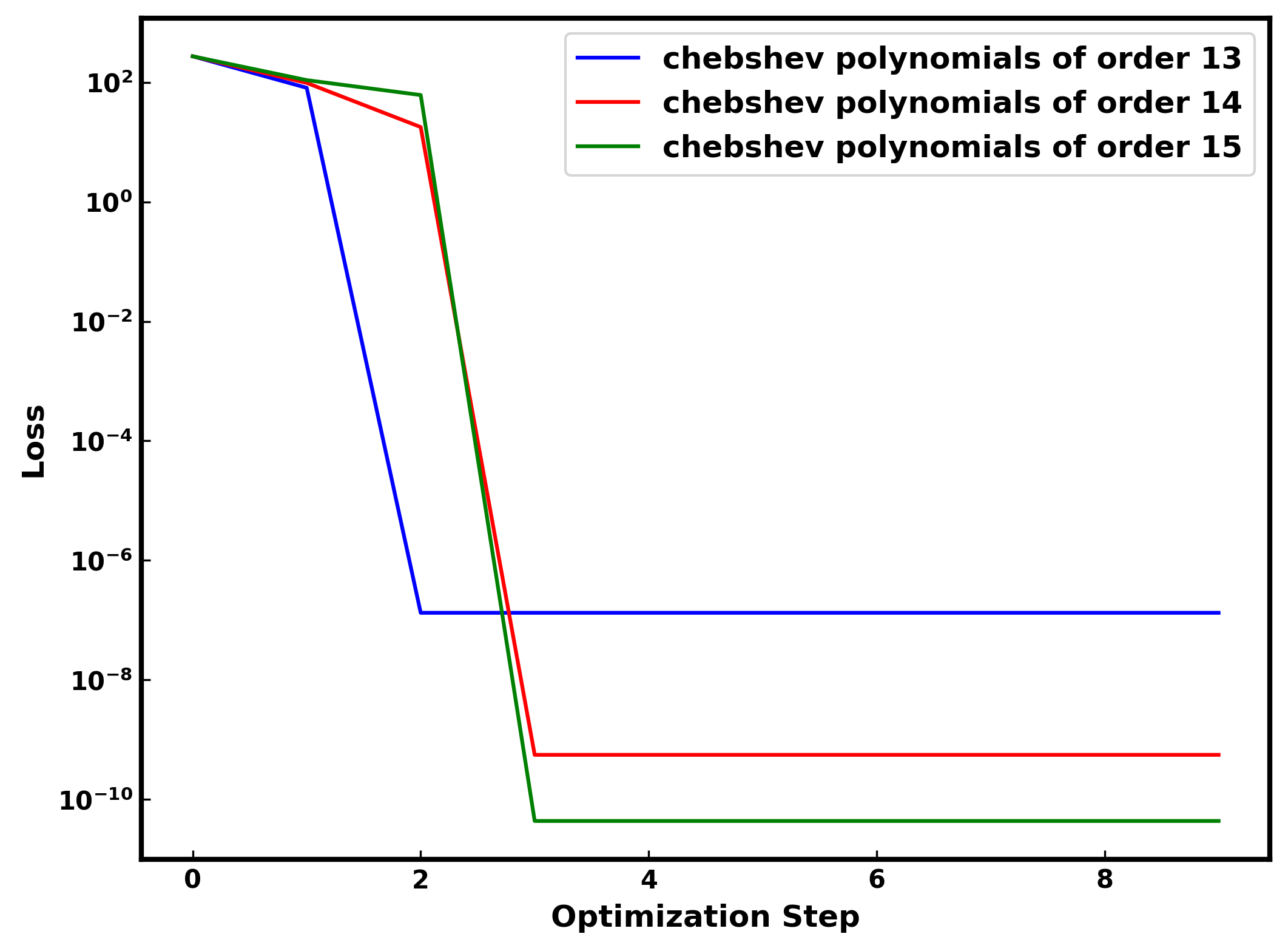}
    \caption{Loss 10}
    \label{loss fig10}
  \end{subfigure}
  \begin{subfigure}[b]{0.30\linewidth}
    \includegraphics[width=\linewidth]{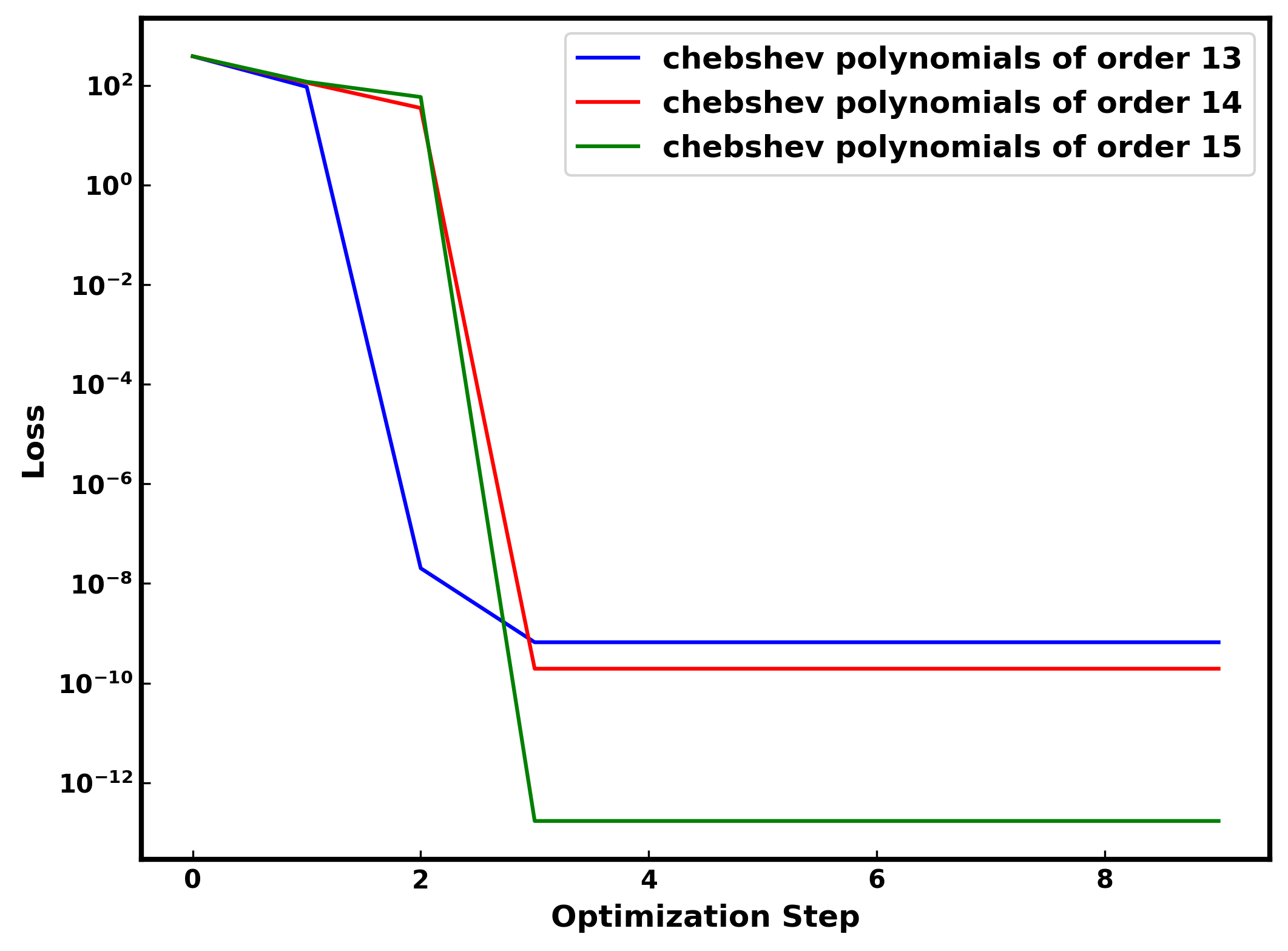}
    \caption{Loss 11}
    \label{loss fig11}
  \end{subfigure}
   \begin{subfigure}[b]{0.30\linewidth}
    \includegraphics[width=\linewidth]{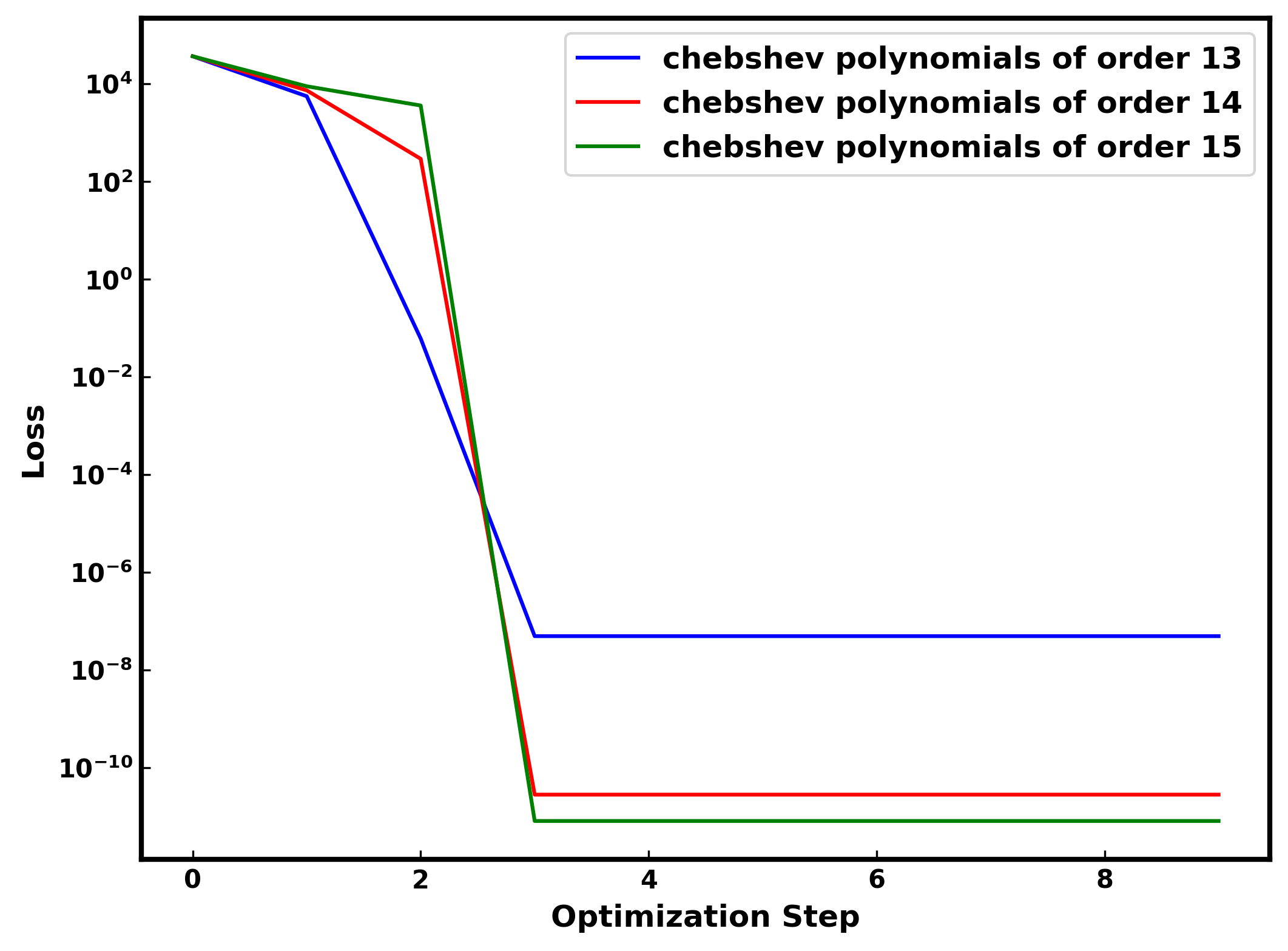}
    \caption{Loss 12}
    \label{loss fig12}
  \end{subfigure}
  \caption{Loss figures of PINN and DFL-TFC frameworks at different parameters for C-S boundary conditions }
   \label{loss figure 2}
\end{figure}

Similarly, details corresponding to Figs \ref{loss fig7} and \ref{loss fig10} are presented in Table \ref{comparison table 4}. Details of Figs \ref{loss fig8} and \ref{loss fig11} are provided in Table \ref{comparison table 5}, whereas those of Figs \ref{loss fig9} and \ref{loss fig12} are summarized in Table \ref{comparison table 6} for C-S boundary conditions.

\begin{table}[H]
\centering
\caption{Comparison of PINN and DFL-TFC for parameters $\alpha = 0.7$, $\gamma = 3$, $\phi = \psi = 0.5$, $N = 1$, $\widetilde{K_p} = 10$, and $\widetilde{q_0} = 2$ (S-S boundary conditions)}
\begin{tabular}{|c c c | c c c|}
\hline
\textbf{PINN} &  &  & \textbf{DFL-TFC} &  &  \\
\hline
\shortstack{No. of \\ Hidden Layers} & Loss & Time & \shortstack{Order of Chebyshev \\ Polynomial} & Loss & Time \\
\hline
1 & \(4.0750\times 10^{-6}\) &  23.86 s & 13 & \(3.1803\times 10^{-8}\) & 4.94 s \\
\hline
2 & \(5.3076\times 10^{-6}\) & 39.43 s & 14 & \(1.8148\times 10^{-10}\) &  5.80 s \\
\hline
3 & \(4.7255\times 10^{-6}\) &  47.43  s & 15 & \(4.8919\times 10^{-10}\) & 7.08 s \\
\hline
\end{tabular}
\label{comparison table 1}
\end{table}

\begin{table}[H]
\centering
\caption{Comparison of PINN and DFL-TFC for parameters $\alpha = 0.8$, $\gamma = 4$, $\phi = \psi = 0.3$, $N = 2$, $\widetilde{K_p} = 5$, and $\widetilde{q_0} = 1$ (S-S boundary conditions)}
\begin{tabular}{|c c c | c c c|}
\hline
\textbf{PINN} &  &  & \textbf{DFL-TFC} &  &  \\
\hline
\shortstack{No. of \\ Hidden Layers} & Loss & Time & \shortstack{Order of Chebyshev \\ Polynomial} & Loss & Time \\
\hline
1 & \( 5.3951\times 10^{-5}\) &  15.60 s & 13 & \(3.4304\times 10^{-10}\) &   5.84 s \\
\hline
2 & \(4.1046\times 10^{-7}\) & 34.86 s & 14 & \(8.7719\times 10^{-11}\) &   6.16 s \\
\hline
3 & \(2.7828\times 10^{-6}\) &  43.26 s & 15 & \(9.9721\times 10^{-11}\) & 8.21 s \\
\hline
\end{tabular}
\label{comparison table 2}
\end{table}

\begin{table}[H]
\centering
\caption{Comparison of PINN and DFL-TFC for parameters $\alpha = 0.2$, $\gamma = 5$, $\phi = \psi = 0.2$, $N = 1$, $\widetilde{K_p} = 8$, and $\widetilde{q_0} = 4$ (S-S boundary conditions)}
\begin{tabular}{|c c c | c c c|}
\hline
\textbf{PINN} &  &  & \textbf{DFL-TFC} &  &  \\
\hline
\shortstack{No. of \\ Hidden Layers} & Loss & Time & \shortstack{Order of Chebyshev \\ Polynomial} & Loss & Time \\
\hline
1 & \(7.9079\times 10^{-5}\) &  36.55 s & 13 & \(1.5806\times 10^{-8}\) &  5.14 s \\
\hline
2 & \(2.2729\times 10^{-3}\) & 57.14 s & 14 & \(3.9759\times 10^{-11}\) &  6.72 s \\
\hline
3 & \(1.4391\times 10^{-4}\) & 58.31 s & 15 & \(6.8713\times 10^{-13}\) & 8.90 s \\
\hline
\end{tabular}
\label{comparison table 3}
\end{table}

\begin{table}[H]
\centering
\caption{Comparison of PINN and Domain mapped FL-TFC for parameters $\alpha = 0.7$, $\gamma = 3$, $\phi = \psi = 0.5$, $N = 1$, $\widetilde{K_p} = 10$, and $\widetilde{q_0} = 2$ (C-S boundary conditions)}
\begin{tabular}{|c c c | c c c|}
\hline
\textbf{PINN} &  &  & \textbf{DFL-TFC} &  &  \\
\hline
\shortstack{No. of \\ Hidden Layers} & Loss & Time & \shortstack{Order of Chebyshev \\ Polynomial} & Loss & Time \\
\hline
1 & \(2.6243\times 10^{-6}\) &  35.78 s & 13 & \(1.3258\times 10^{-7}\) & 4.80 s \\
\hline
2 & \(5.8697\times 10^{-5}\) &  56.07 s & 14 & \(5.5581\times 10^{-10}\) &   7.11 s \\
\hline
3 & \(1.2247\times 10^{-6}\) &  49.03 s & 15 & \(4.3449\times 10^{-11}\) &9.44 s \\
\hline
\end{tabular}
\label{comparison table 4}
\end{table}

\begin{table}[H]
\centering
\caption{Comparison of PINN and Domain mapped FL-TFC for parameters $\alpha = 0.8$, $\gamma = 4$, $\phi = \psi = 0.3$, $N = 2$, $\widetilde{K_p} = 5$, and $\widetilde{q_0} = 1$ (C-S boundary conditions)}
\begin{tabular}{|c c c | c c c|}
\hline
\textbf{PINN} &  &  & \textbf{DFL-TFC} &  &  \\
\hline
\shortstack{No. of \\ Hidden Layers} & Loss & Time & \shortstack{Order of Chebyshev \\ Polynomial} & Loss & Time \\
\hline
1 & \(1.6465\times 10^{-6}\) &  27.36 s & 13 & \(6.6098\times 10^{-10}\) &  6.39 s \\
\hline
2 & \(2.7963\times 10^{-5}\) &  45.26 s & 14 & \(1.9484\times 10^{-10}\) &  8.09 s \\
\hline
3 & \(1.1195\times 10^{-5}\) & 41.58 s & 15 & \(1.7089\times 10^{-13}\) & 9.66 s \\
\hline
\end{tabular}
\label{comparison table 5}
\end{table}

\begin{table}[H]
\centering
\caption{Comparison of PINN and Domain mapped FL-TFC for parameters $\alpha = 0.2$, $\gamma = 5$, $\phi = \psi = 0.2$, $N = 1$, $\widetilde{K_p} = 8$, and $\widetilde{q_0} = 4$ (C-S boundary conditions)}
\begin{tabular}{|c c c | c c c|}
\hline
\textbf{PINN} &  &  & \textbf{DFL-TFC} &  &  \\
\hline
\shortstack{No. of \\ Hidden Layers} & Loss & Time & \shortstack{Order of Chebyshev \\ Polynomial} & Loss & Time \\
\hline
1 & \(8.8678\times 10^{-4}\) &  35.76 s & 13 & \(4.9020\times 10^{-8}\) & 6.71 s \\
\hline
2 & \(2.7195\times 10^{-4}\) &  58.58 s & 14 & \(2.7857\times 10^{-11}\) & 8.23 s \\
\hline
3 & \(3.7222\times 10^{-1}\) &  71.00 s & 15 & \(8.0515\times 10^{-12}\) & 10.34 s \\
\hline
\end{tabular}
\label{comparison table 6}
\end{table}

\section{Conclusion}\label{sec. conclusion}

This work introduces a robust DFL-TFC based framework for analyzing bending behavior of a tapered perforated beam subjected to an exponential load under simply supported (S-S) and clamped simply supported (C-S) boundary conditions.

The DFL-TFC method demonstrates excellent computational performance for finding bending deflection. This approach ensures high accuracy while achieving faster convergence and reduced computational cost compared to PINN method. Moreover, the inherent satisfaction of boundary conditions within the DFL-TFC framework enhances stability and reliability of the obtained solutions.
Overall, the proposed framework proves to be an efficient and powerful tool for analyzing complex perforated beam structures. Methodology can be readily extended to more advanced structural configurations and loading conditions, offering significant potential for future research in computational mechanics and engineering applications.

New outcomes reveal that the bending deflection decreases with an increase $\alpha$, as a higher material presence enhances stiffness of beam. Similarly, increasing $\phi$ and $\psi$ leads to reduced deflection due to the improved stiffness resulting from the modified geometry and mass distribution. In contrast, deflection increases with higher values of  $\gamma$, as it intensifies the applied load. Furthermore, presence of an elastic foundation plays a significant role, where an increase in $\widetilde{K_p}$ provides greater resistance to deformation, thereby reducing the overall deflection.

\subsection*{Acknowledgment}
The first author expresses gratitude to Council of Scientific and Industrial Research, India for providing fellowship.

\subsection*{Appendix} \label{appendix}
Here, details of CE are presented for both S-S and C-S boundary conditions. The CE is expressed as
\begin{equation}
    \widehat{\overline{W}_p}(X) = h(X) + \mathcal{M}(\boldsymbol{\rho}) \Psi_X,
\end{equation}

where 
\begin{equation}
\Psi_X = \left[ 1 \hspace{0.3cm} {}^X\psi_1 \hspace{0.3cm} {}^X\psi_2 \hspace{0.3cm} {}^X\psi_3 \hspace{0.3cm} {}^X\psi_4\right]^T
\end{equation}
is vector of switching functions \({}^X\psi_1\), \({}^X\psi_2\), \({}^X\psi_3\) and \({}^X\psi_4\) corresponds to the constraints at \(X = 0\) and \(X = 1\).\\
Support functions are choosen as 
\begin{equation}
s_1(X) = 1, \quad
s_2(X) = X, \quad
s_3(X) = X^2, \quad
s_4(X) = X^3
\end{equation}

\begin{equation}
\left[\alpha_{ij}\right]
=
\left[{}^X\mathfrak{C}^i[s_j(X)]\right]^{-1}
=
\begin{bmatrix}
\mathtt{s}_1(0) & \mathtt{s}_2(0) & \mathtt{s}_3(0) & \mathtt{s}_4(0) \\[6pt]
\mathtt{s}_1(1) & \mathtt{s}_2(1) & \mathtt{s}_3(1) & \mathtt{s}_4(1) \\[6pt]
\mathtt{s''}_1(0) & \mathtt{s''}_2(0) & \mathtt{s''}_3(0) & \mathtt{s''}_4(0) \\[6pt]
\mathtt{s''}_1(1) & \mathtt{s''}_2(1) & \mathtt{s''}_3(1) & \mathtt{s''}_4(1)
\end{bmatrix}^{-1}
\end{equation}

\begin{equation}
=
\begin{bmatrix}
1 & 0 & 0 & 0\\
1 & 1 & 1 & 1\\
0 & 0 & 2 & 0 \\
0 & 0 & 2 & 6
\end{bmatrix}^{-1}
\end{equation}

\begin{equation}    
=
\begin{bmatrix}
1 & 0 & 0 & 0\\
-1 & 1 & -\frac{1}{3} & -\frac{1}{6}\\
0 & 0 & \frac{1}{2} & 0 \\
0 & 0 & -\frac{1}{6} & \frac{1}{6}
\end{bmatrix}
\end{equation}

Switching functions in \(X\) are given by     
\begin{equation*}
{}^{X}\psi_1 = \mathtt{s}_1(X)\cdot 1 +\mathtt{s}_2(X)\cdot (-1) + \mathtt{s}_3(X)\cdot 0 + \mathtt{s}_4(X)\cdot 0  = 1 - X 
\end{equation*}

\begin{equation*}
 {}^{X}\psi_2  = \mathtt{s}_1(X)\cdot 0 + \mathtt{s}_2(X)\cdot 1 + \mathtt{s}_3(X)\cdot 0 + \mathtt{s}_4(X)\cdot 0 = X   
\end{equation*}

\begin{equation*}
  {}^{X}\psi_3
= \mathtt{s}_1(X)\cdot 0
+ \mathtt{s}_2(X)\cdot \left(-\frac{1}{3}\right)
+ \mathtt{s}_3(X)\cdot \frac{1}{2}
+ \mathtt{s}_4(X)\cdot \left(-\frac{1}{6}\right)
= -\frac{1}{3}X + \frac{1}{2}X^2 - \frac{1}{6}X^3  
\end{equation*}
\begin{equation*}
  {}^{X}\psi_4
= \mathtt{s}_1(X)\cdot 0
+ \mathtt{s}_2(X)\cdot \left(-\frac{1}{6}\right)
+ \mathtt{s}_3(X)\cdot 0
+ \mathtt{s}_4(X)\cdot \left(\frac{1}{6}\right)
= -\frac{1}{6}X + \frac{1}{6}X^3  
\end{equation*}

CE for the S-S boundary conditions is obtained as follows: 
\begin{equation}
\boxed{
\widehat{\overline{W_p}}(X) = h(X) + (X - 1) h(0) - X h(1) + \frac{2X - 3X^2 + X^3}{6}h''(0)+ \frac{X - X^3}{6}h''(1).}
\end{equation}
Similarly, CE for the C-S boundary conditions is given below: 

 \begin{equation}
\boxed{
\begin{aligned}  
 \widehat{\overline{W_p}}(X) = h(X) + \frac{-2 + 3X^2 - X^3}{2}h(0)+ \frac{-2X + 3X^2 - X^3}{2}h'(0) \\
 + \frac{X^3 - 3X^2}{2}h(1) + \frac{X^2 - X^3}{4}h''(1).
 \end{aligned}
 }
\end{equation}

All the experiments are conducted on the Google Colab platform using Python with CPU-based hardware.

\bibliographystyle{ieeetr}
\bibliography{col_reference}
\end{document}